\documentclass[a4paper,12pt]{article}
   \usepackage{hyperref}
   \usepackage{amsmath}
   \usepackage{amsfonts}
   \usepackage{amssymb}
   \usepackage{mathrsfs}
   \usepackage{graphicx}
   \usepackage[abs]{overpic}
   \usepackage{float}

\newtheorem{thm}{Theora}[section]
\newtheorem{theo}[thm]{Theorem}

\newtheorem{coro}[thm]{Corollary}
\newtheorem{lemm}[thm]{Lemma}
\newtheorem{prop}[thm]{Proposition}

\newtheorem{rem}[thm]{Remark}

\newenvironment{Proof}{\smallskip\noindent{\bf Proof.}\rm}
{\hfill $\Box$\medskip}
\newenvironment{proof}{\smallskip\noindent{\bf Proof}\rm}
{\hfill $\Box$\medskip}
\renewcommand\({\left(}
\renewcommand\){\right)}
\renewcommand\[{\left[}
\renewcommand\]{\right]}

\newcommand\la{\lambda}
\newcommand\ga{\gamma}
\newcommand{\be}{\begin{equation}}
\newcommand{\ee}{\end{equation}}
\newcommand{\ba}{\begin{array}}
\newcommand{\ea}{\end{array}}
\newcommand{\bea}{\begin{eqnarray*}}
\newcommand{\eea}{\end{eqnarray*}}
\newcommand{\bean}{\begin{eqnarray}}
\newcommand{\eean}{\end{eqnarray}}
\newcommand\ep{\varepsilon}

\newcommand\D{Dim}
\makeatletter \@addtoreset{equation}{section}

\makeatother
\begin{document}
%\Large\bf{
\title{Boundary Controllability Of Two Vibrating Strings Connected By A Point Mass With Variable Coefficients}
%Exact controllability of a one-dimensional wave equation with an
%interior point mass and variable coefficients.
\author{Jamel Ben Amara \thanks{Department of Mathematics, Faculty of Sciences of Tunis, University of Tunis el Manar,
Mathematical Engineering Laboratory, Tunisia; e-mail:
jamel.benamara@fsb.rnu.tn.}~~~~~~~~~~Emna Beldi
\thanks{Department of Mathematics, Tunisia Polytechnic School, Mathematical Engineering Laboratory,
University of Carthage, Tunisia; e-mail: em.beldi@gmail.com.}}
\date{}
\maketitle {\bf Abstract.} S. Hansen and E. Zuazua [SIAM J. Cont.
Optim., 1995] studied the problem of exact controllability of two
strings connected by a point mass with constant physical
coefficients. In this paper we study the same problem with variable
physical coefficients. This system is generated by the following
equations
$$\rho(x) u_{tt}=(\sigma(x) u_{x})_{x}-q(x)u,~~~~x\in (-1,0)\cup (0,1),~t>0,$$
$$Mu_{tt}(0,t)+\sigma_{1}(0)u_{x}(0^{-},t)-\sigma_{2}(0)u_{x}(0^{+},t)=0,~~~t>0,$$
with Dirichlet boundary condition on the left end and a control acts
on the right end. We prove that this system is exactly controllable
in an asymmetric space for the control time $T>
2\int_{-1}^{1}(\frac{\rho(x)}{\sigma(x)})^{\frac{1}{2}}dx$. We
establish the equivalence between a suitable asymmetric norm of the
initial data and the $L^{2}(0,T)$-norm of $u_{x}(1,t)$ (where $u$ is
the solution of the uncontrolled system). Our approach is mainly
based on a detailed spectral analysis and the theory of divided
differences. More precisely, we prove that the
spectral gap tends to zero with a precise asymptotic estimate.\\
{\bf Key words.} Boundary control, point mass, Riesz basis,
vibrating string, variable coefficients.\\
{\bf Mathematics Subject Classification.} 35P, 47A, 93B.

%%%%%%%%%%%%%%%%%%%%%%%%%%%%%%%%%%%%%%%%%%%%%%%%%%%%%%%%%%%%%%%%%%%%%%%%%%%%%%%%%%%%%%%%%%%%%%%%%%%%%
%%%%%%%%%%%%%%%%%%%%%%%%%%%%%%%%%%%%%%%%%%%%%%%%%%%%%%%%%%%%%%%%%%%%%%%%%%%%%%%%%%%%%%%%%%%%%%%%%%%%%
\section{Introduction} The controllability of mechanical structures with attached
masses has been extensively investigated for several decades. To
mention some examples, see \cite{han} for heat equation connected by
a point mass, \cite{CAS, E.Z, Rouch} for strings with an interior
point mass  and \cite{Beam2, Beam1, Beam3} for beams with an
attached point mass. See also \cite{ee, ren1, ren2} for
networks of strings or beams with attached masses.\\
In this paper we study the boundary controllability of a linear
hybrid system modeling two vibrating non-homogeneous strings
connected by a point mass. By means of the functions:
$$u=u(x,t), x \in (-1,0), t>0,$$
$$v=v(x,t), x\in (0,1), t>0,$$
we describe the vertical displacements of the first and the second
string, respectively. The position of the mass $M>0$ attached to the
strings at the point $x=0$ is described by the function $z=z(t)$ for $t>0$.\\
The linear system modeling the vibrations of these strings is given
as follows
\begin{equation}\label{a}
\begin{cases}
\rho_{1}(x)u_{tt}=(\sigma_{1}(x)u_{x})_{x}-q_{1}(x)u,~~~~~~~~~~~~~~~~~&x\in (-1, 0), t>0,\\
\rho_{2}(x)v_{tt}=(\sigma_{2}(x)v_{x})_{x}-q_{2}(x)v,~~~~~~~~~~~~~&x\in (0, 1), t>0, \\
u(0,t)=v(0,t)=z(t),~~~~~~~~~~~~~~~~~~~~~~~~~~~~~~~~~~&t>0, \\
Mz_{tt}(t)+\sigma_{1}(0)u_{x}(0,t)-\sigma_{2}(0)v_{x}(0,t)=0, ~~~~~~~~~&t>0, \\
\end{cases}
\end{equation}
with the following boundary conditions
\begin{equation}\label{22}
u(-1,t)=v(1,t)=0,~for~t>0
\end{equation}
and the following initial conditions
\begin{equation}\label{j}
\begin{cases}
  u(x,0)=u^{0}(x), u_{t}(x,0)=u^{1}(x),~~~~~~~~~~~~~~~~~~~~~~&x\in(-1,0), \\
  v(x,0)=v^{0}(x), v_{t}(x,0)=v^{1}(x),~~~~~~~~~~~~~~~~~~~~~~~&x\in(0,1), \\
  z(0)=z^{0}, z_{t}(0)=z^{1}. \\
\end{cases}
\end{equation}
The coefficients $\rho_{i}(x)$ and $\sigma_{i}(x)$ represent
respectively the density and the tension of each string, $i=1,2$.
The potentials are denoted by the functions $q_{1}(x)$ and
$q_{2}(x)$. In this study we are interested in the boundary
controllability when a control function $p=p(t)$ acts on the system
through the extreme $x=1$. Then the boundary conditions \eqref{22}
is replaced by
\begin{equation}\label{b1}
u(-1,t)=0,~ v(1,t)=p(t),~~t>0.
\end{equation}
Throughout this paper the coefficients $\rho_{i}$ and $\sigma_{i}$
are assumed to be uniformly positive, $q_{i}$ is nonnegative (i=1,2)
and
$$\rho_{1},~\sigma_{1} \in H^{2}(-1,0),~q_{1}\in L^{1}(-1,0),$$
$$\rho_{2},~\sigma_{2} \in H^{2}(0,1),~q_{2}\in L^{1}(0,1).$$
In what follows, we introduce the following spaces
\begin{equation}\label{h}
\mathcal{H}_{0}=L^{2}(-1,0)\times L^{2}(0,1)\times \mathbb{R}.
\end{equation}
\bean
&&\mathcal{V}_{1}=\{u\in H^{1}(-1,0)~~|~~ u(-1)=0\}\nonumber\\
&&\mathcal{V}_{2}=\{v\in H^{1}(0,1)~~|~~ v(1)=0\}\label{w3}\\
&&\mathcal{V}=\{(u,v)\in \mathcal{V}_{1}\times
\mathcal{V}_{2}~~|~~u(0)=v(0)\}\nonumber\eean and
\begin{equation}\label{h1}
\mathcal{H}_{-1}=\mathcal{V}_{1}'\times \mathcal{V}_{2}'\times
\mathbb{R},
\end{equation}
where the space $(\mathcal{V}_{i})'$ denotes the dual space of
$\mathcal{V}_{i}$, $(i=1,2)$.\\
In a similar way as in \cite{E.Z}, we have the following existence
and uniqueness result for System \eqref{a}-\eqref{j}-\eqref{b1}:
\begin{prop} For every $p \in L^{2}(0,T)$,
$(u^{0},v^{0},z^{0}) \in \mathcal{H}_{0}$ and $(u^{1},v^{1},z^{1})
\in \mathcal{H}_{-1}$ such that $(u^{0},u^{1})\in
\mathcal{V}_{1}\times L^{2}(-1,0)$ and $u^{0}(0)=z^{0}$, there
exists a solution of System \eqref{a}-\eqref{b1} with initial data
\eqref{j} in the class:
$$(u,v,z) \in C([0,T],\mathcal{H}_{0})\cap C^{1}([0,T],\mathcal{H}_{-1})$$ and
$$u\in C([0,T],\mathcal{V}_{1})\cap C^{1}([0,T],L^{2}(-1,0)).$$
\end{prop}
Now, we can state our main result in this paper:
\begin{theo}\label{hj}
Suppose $$T>2(\gamma_{1}+\gamma_{2})
=2\left(\int_{-1}^{0}\left(\frac{\rho_{1}(x)}{\sigma_{1}(x)}\right)^{\frac{1}{2}}dx+
\int_{0}^{1}\left(\frac{\rho_{2}(x)}{\sigma_{2}(x)}\right)^{\frac{1}{2}}dx\right).$$
Then, for every $(u^{0},v^{0},z^{0}) \in \mathcal{H}_{0}$ and
$(u^{1},v^{1},z^{1}) \in \mathcal{H}_{-1}$ such that\\
$(u^{0},u^{1})\in \mathcal{V}_{1}\times L^{2}(-1,0)$ and
$u^{0}(0)=z^{0}$, there exists a control $p \in L^{2}(0,T)$ such
that the solution of \eqref{a}-\eqref{j}-\eqref{b1} satisfies
\begin{equation*}
\left\{
\begin{array}{l}
  u(x,T)=u_{t}(x,T)=0, \\
  v(x,T)=v_{t}(x,T)=0, \\
  z(T)=z_{t}(T)=0.
\end{array}
\right.
\end{equation*}
\end{theo}
For the proof of this theorem, we establish the following result: if
$T>2(\gamma_{1}+\gamma_{2})$, then for every solution $(u,v,z)$ of
the uncontrolled problem \eqref{a}-\eqref{22}-\eqref{j} there exist
two positive constants $C_{1}$ and $C_{2}$ such that
\begin{equation}\label{obs}
C_{1} \|U^{0}\|^{2}_{\mathcal{Y}}\leq
\|v_{x}(1,t)\|^{2}_{L^{2}(0,T)}\leq C_{2}
\|U^{0}\|^{2}_{\mathcal{Y}},
\end{equation}
where $U^{0}=(u^{0},v^{0},z^{0}),(u^{1},v^{1},z^{1})$ and
\begin{align}
\mathcal{Y}=&\big\{U^{0}=((u^{0},v^{0},z^{0}),(u^{1},v^{1},z^{1}))\in
\mathcal{H}_{0}\times \mathcal{H}_{-1}|~(v^{0},v^{1})\in
\mathcal{V}_{2}\times L^{2}(0,1)\big\}.\label{Y}
\end{align}
In \cite{E.Z}, S. Hansen and E. Zuazua studied System
\eqref{a}-\eqref{22}, for constant physical parameters $\rho_{i}>0$,
$\sigma_{i}>0$ and $q_{i}\equiv0$, $(i=1,2)$ when the two strings
occupy the intervals $(-\ell_{1},0)$ and $(0,\ell_{2})$
$(\ell_{1},\ell_{2}>0)$. Using the Hilbert Uniqueness Method "HUM"
(cf. J.L. Lions \cite{J.L2}), they proved the exact controllability
in an asymmetric space when the control acts at one end of the
string-mass-system for a control time
$T>2(\ell_{1}\sqrt{\frac{\rho_{1}}{\sigma_{1}}}+
\ell_{2}\sqrt{\frac{\rho_{2}}{\sigma_{2}}})$. In the particular
case, when $\rho_{i}=\sigma_{i}=\ell_{i}=1$ and $q_{i}\equiv0$
(i=1,2), the proof of the exact controllability is mainly based on
the theory of non-harmonic Fourier series. Since in this case the
spectral gap tends to zero, they employed a result of D. Ullrich
\cite{D.U} to prove the observability inequality. At the end of
their paper, they considered the case of variable coefficients
$\rho_{i}$, $\sigma_{i}$ and $q_{i}$ $(i=1,2)$, where Theorem
\ref{hj} was enunciated in a form of a conjecture. Later C. Castro
in \cite{CAS} proved the same result by using a different approach,
he showed that the solutions of the above system can be viewed as
weak limits of solutions of a string equation with unit density on
$(-1 ,1)\setminus (-\epsilon, \epsilon)$ and with density
$\frac{1}{2\epsilon}$ on $(-\epsilon, \epsilon)$. Recently, S.
Avdonin and J. Edward in \cite{Avd} have studied the exact
controllability of a vibrating string with $N$ attached masses
inside an interval $(0,\ell)$ with a Dirichlet control at one end.
Their problem in the case of a string with a single attached mass,
coincides with System \eqref{a}-\eqref{22} for
$\rho_{i}=\sigma_{i}=1$ (i=1,2). Using a different approach that
used in \cite{CAS, E.Z}, they proved the exact controllability in an
asymmetric space for $T\geq 2\ell$. In particular, they established
the equivalence between an asymmetric norm of the
initial data and the $L^{2}(0,T)$-norm of the control $p(t)$.\\
Our approach is essentially based on a precise computation of the
associated spectral gap together with a suitable non-harmonic
Fourier series result for the time exponentials. We show that there
exists a subsequence of eigenfrequencies
$(\sqrt{\lambda_{\varphi(n)}})$ associated with System
\eqref{a}-\eqref{22} such that
$\sqrt{\lambda_{\varphi(n)+1}}-\sqrt{\lambda_{\varphi(n)}}
=\mathcal{O}(\frac{1}{n})$
and\\$\sqrt{\la_{n+2}}-\sqrt{\la_{n}}>2\delta, ~n\geq1$, for some
$\delta>0$. In the process of the last result, we establish an
interpolation formula between the eigenvalues of System
\eqref{a}-\eqref{22} and those of the regular problem
\eqref{a}-\eqref{22} for $M=0$. Instead of Ullrich's theorem used in
\cite{CAS, E.Z}, we here apply a result of C. Baiocchi et al.
\cite{Bi1}. More general results concerning divided differences
theory can be
found in \cite{AVT} and \cite{Ah}.\\ %In Remark \ref{rem}, we noted
%that in the case of symmetric coefficients the eigenvalues of the
%Dirichlet boundary problems defined in each interval $[-1,0]$ and
%$[0,1]$ coincide. Hence, it is possible to use Ullrich's
%theorem as in \cite{CAS} and \cite{E.Z}.\\
The rest of the paper is divided in the following way: In section 2
we associate to System \eqref{a}-\eqref{22} a self-adjoint operator
defined in the Hilbert space $\mathcal{H}_{0}$ (defined by
\eqref{h}. Moreover, we give a result on the well-posedness of
system \eqref{a}-\eqref{22}-\eqref{j}. In the next section we
establish the asymptotic properties of the associated spectral gap
and the asymptotes of the eigenfunctions. In section 4 we define an
asymmetric space and we give its characterization. In the last
section we prove Theorem \ref{hj} and the observability inequality
\eqref{obs}.
\section{Operator Framework And Well-posedness}
We begin with considering the following closed subspace of
$\mathcal{V} \times \mathbb{R}$, where $\mathcal{V}$ is defined by
\eqref{w3}
\begin{equation}\label{v}
\mathcal{W}=\{(u,v,z)\in \mathcal{V} \times \mathbb{R}:
u(0)=v(0)=z\}
\end{equation}
equipped with the norm
\begin{equation}\label{j33}
\|(u,v,z)\|^{2}_{\mathcal{W}}=\|u\|^{2}_{\mathcal{V}_{1}}+
\|v\|^{2}_{\mathcal{V}_{2}},
\end{equation}
where
\begin{equation}\label{j11}
\|u\|^{2}_{\mathcal{V}_{1}}=\int_{-1}^{0}|u_{x}(x)|^{2}dx,~~
\|v\|^{2}_{\mathcal{V}_{2}}=\int_{0}^{1}|v_{x}(x)|^{2}dx\end{equation}
%\begin{equation}\label{j33}
%\|(u,v)\|^{2}_{\mathcal{V}}=\|u\|^{2}_{\mathcal{V}_{1}}+
%\|v\|^{2}_{\mathcal{V}_{2}},\end{equation}
and the spaces $\mathcal{V}_{1}$ and $\mathcal{V}_{2}$ are
defined by \eqref{w3}.\\
In the sequel we introduce the operator $\mathcal{A}$ defined in
$\mathcal{H}_{0}$, which is defined by \eqref{h}, by setting:
\begin{equation}\label{est}
\mathcal{A}Y=\begin{cases}
\frac{1}{\rho_{1}}(-(\sigma_{1}u')'+q_{1}u),\\
\frac{1}{\rho_{2}}(-(\sigma_{2}v')'+q_{2}v),\\
\frac{1}{M}(\sigma_{1}(0)u'(0)-\sigma_{2}(0)v'(0)),
\end{cases}
\end{equation}
where $Y=(u,v,z)^{t}$ on the domain
\begin{equation*}
D(\mathcal{A})=\{Y|Y=(u,v,z)\in \mathcal{W}:~(u,v)\in
H^{2}(-1,0)\times H^{2}(0,1)\},
\end{equation*}
which is dense in $\mathcal{H}_{0}$. Note that the space
$\mathcal{H}_{0}$ is equipped with the scalar product
$\langle./.\rangle_{\mathcal{H}_{0}}$ defined by: for all
$Y_{1}=(u_{1},v_{1},\alpha_{1})^{t}$ and
$Y_{2}=(u_{2},v_{2},\alpha_{2})^{t}\in \mathcal{H}_{0}$, where
$^{t}$ denotes the transposition, we have
\begin{equation*}
\langle Y_{1},
Y_{2}\rangle_{\mathcal{H}_{0}}=\int_{-1}^{0}u_{1}(x)u_{2}(x)\rho_{1}(x)dx+
\int_{0}^{1}v_{1}(x)v_{2}(x)\rho_{2}(x)dx+M\alpha_{1}\alpha_{2}.
\end{equation*}
\begin{prop}\label{rr}
The linear operator $\mathcal{A}$ is self-adjoint and positive such
that $\mathcal{A}^{-1}$ is compact. Moreover, the operator
$\mathcal{A}^{\frac{1}{2}}$ generates a strongly continuous
semi-group on $\mathcal{H}_{0}$.
\end{prop}
\begin{Proof}
Let $Y=(u,v,z)^{t}\in D(\mathcal{A})$, then by a simple integration
by parts, we have
\begin{align*}
\langle \mathcal{A}Y, Y\rangle_{\mathcal{H}_{0}}
=\int_{-1}^{0}\sigma_{1}(x)|u^{'}|^{2}dx+\int_{-1}^{0}q_{1}(x)|u|^{2}dx
+\int_{0}^{1}\sigma_{2}(x)|v^{'}|^{2}dx+\int_{0}^{1}q_{2}(x)|v|^{2}dx.
\end{align*}
Since $\sigma_{i}>0$ and $q_{i}>0$, (i=1,2), then $\langle
\mathcal{A}Y, Y\rangle_{\mathcal{H}_{0}}>0$, and hence, the operator
$\mathcal{A}$ is symmetric on $\mathcal{H}_{0}$. To prove that this
operator is self-adjoint, it suffices to show that
$Ran(\mathcal{A}-iId)=\mathcal{H}_{0}$, and this is omitted here.\\
It is easy to see that the space $\mathcal{W}$ is continuously and
compactly embedded in the space $\mathcal{H}_{0}$, and hence, the
operator $\mathcal{A}^{-1}$ is compact in $\mathcal{H}_{0}$.
\end{Proof}\\
According to Proposition \ref{rr}, the operator $\mathcal{A}$ is
positive and self-adjoint, and hence, it generates a scale of
Hilbert spaces $\mathcal{H}_{\theta}=D(\mathcal{A}^{\theta})$,
$\theta \in \mathbb{R}$. In particular, $\mathcal{H}_{0}$ is defined
by \eqref{h} and
$\mathcal{H}_{\frac{1}{2}}=D(\mathcal{A}^{\frac{1}{2}})$. The norm
in $\mathcal{H}_{\frac{1}{2}}$ is given by
\begin{align*}
\|y\|_{\frac{1}{2}}^{2}&=\langle\mathcal{A}^{\frac{1}{2}}y,
\mathcal{A}^{\frac{1}{2}}y\rangle_{\mathcal{H}_{0}},\\
&=\int_{-1}^{0}(\sigma_{1}(x)|u'|^{2}+q_{1}(x)|u|^{2})dx
+\int_{0}^{1}(\sigma_{2}(x)|v'|^{2}+q_{2}(x)|v|^{2})dx.
\end{align*}
Our assumptions on $\sigma_{i}(x)$ and $q_{i}(x)$ $(i=1,2)$ imply
that $\mathcal{H}_{\frac{1}{2}}$ is topologically equivalent to the
subspace $\mathcal{W}$.
\begin{prop}\label{33}
The map $Y=(u,v,z)\rightarrow (u,v)$ is an homeomorphism from
$D(\mathcal{A}^{\frac{1}{2}})$ into $\mathcal{V}$, where
$\mathcal{V}$ is defined by \eqref{w3}.
\end{prop}
 %\begin{Proof} A similar result is proved at the end of the paper
% \cite{ab}.
%\end{Proof}\\
%\subsection{Homogeneous boundary conditions}
Obviously, the Cauchy problem \eqref{a}-\eqref{j} can be rewritten
in the abstract form
\begin{equation*}
\ddot{Y}=-\mathcal{A}Y,~(Y(0),Y_{t}(0))=(Y^{0},Y^{1}),
\end{equation*}
where $\mathcal{A}$ is defined by \eqref{est},
$Y^{0}=(u^{0},v^{0},z^{0})^{t}$ and $Y^{1}=(u^{1},v^{1},z^{1})^{t}$.
As a consequence of Proposition \ref{rr}, we have the following
existence and uniqueness result for Problem
\eqref{a}-\eqref{22}-\eqref{j}:
\begin{theo}
For every $((u^{0},v^{0}),(u^{1},v^{1},z^{1})) \in \mathcal{V}
\times \mathcal{H}_{0}$, there exists a unique solution of System
\eqref{a}-\eqref{22} with initial conditions \eqref{j} in the class
$$(u,v,z)\in C([0,T],\mathcal{W})\cap C^{1}([0,T], \mathcal{H}_{0}).$$
%Or equivalently, for every $((u^{0},v^{0}),(u^{1},v^{1},z^{1})) \in
%\mathcal{V} \times \mathcal{H}$, Problem \eqref{jik} have a unique
%solution.
\end{theo}
\section{Spectral Gap And Asymptotic Of The Eigenfunctions} In
this section we investigate the spectral properties of the operator
$\mathcal{A}$, in particular we establish the asymptotic behavior of
the spectral gap $\sqrt{\la_{n+1}}-\sqrt{\la_{n}}$. Similar
asymptotic estimate in the case of constant coefficients was
obtained in \cite{CAS, E.Z}. We consider the following spectral
problem which arises by applying separation of variables to System
\eqref{a}-\eqref{22}:
\begin{equation}\label{an}
\begin{cases}
-(\sigma_{1}(x)u'(x))'+q_{1}(x)u(x)=\lambda \rho_{1}(x)u(x),~~~~~~~~~~~~&x\in (-1,0), \\
  -(\sigma_{2}(x)v'(x))'+q_{2}(x)v(x)=\lambda \rho_{2}(x)v(x),~~~~~~~~~~~~~~~~&x\in(0,1),\\
  u(-1)=v(1)=0,~u(0)=v(0),\\
  \sigma_{1}(0)u'(0)-\sigma_{2}(0)v'(0)=\lambda M u(0).
\end{cases}
\end{equation}
\begin{lemm}\label{lemm2}
The spectrum of the operator $\mathcal{A}$ is discrete. It consists
of an increasing sequence of positive and simple eigenvalues
$(\lambda_{n})_{n\in\mathbb{N}^{*}}$ tending to $+\infty$:
$$0<\lambda_{1}<\lambda_{2}<.......<\lambda_{n}<.....\underset{n\rightarrow +\infty}{\longrightarrow}+\infty.$$
Moreover, the corresponding eigenfunctions
$\Phi_{n}(x)=(\phi_{n}(x),\phi_{n}(0))$ form an orthogonal basis in
$\mathcal{H}_{0}$, with $\phi_{n}$ are the eigenfunctions of the
eigenvalue problem \eqref{an}.
\end{lemm}
\begin{Proof}
Here we have only to prove the simplicity of the eigenvalues
$\la_{n}$ for all $n\in\mathbb{N^{*}}$. Let $\tilde{u}(x,\la)$ and
$\tilde{v}(x,\la)$ be the solutions of the initial value problems
\begin{equation}\label{w1}
\left\{
\begin{array}{lll}
-(\sigma_{1}(x)u')'+q_{1}(x)u=\lambda \rho_{1}(x) u, ~~x\in (-1,\,0),\\
u(-1)=0,~u'(-1)=1
\end{array}
\right.
\end{equation}
and
\begin{equation}\label{w2}
\left\{
\begin{array}{lll}
-(\sigma_{2}(x)v')'+q_{2}(x)v=\lambda \rho_{2}(x)v,~~x\in (0,\,1),\\
v(1)=0,~v'(1)=-1,
\end{array}
\right.
\end{equation}
respectively. Let $\lambda$ be an eigenvalue of the operator
$\mathcal{A}$ and $E_{\lambda}$ be the corresponding eigenspace. For
any eigenfunction $\phi (x,\lambda)$ of $E_{\lambda}$, $\phi
(x,\lambda)$ can be written in the form
\begin{equation}\label{8} \phi(x,\la)= \left\{
\begin{array}{lll}
c_{1}\tilde{u}(x,\la),&-1\leq x\leq 0,\\
c_{2}\tilde{v}(x,\la),&0\leq x\leq 1,\\
\end{array}
\right.
\end{equation}
where $c_{1}$ and $c_{2}$ are two constants. The first condition at
$x=0$ of Problem \eqref{an} is equivalent to
$$c_{1}\tilde{u}(0,\la)=c_{2}\tilde{v}(0,\la).$$
If $\tilde{u}(0,\la)\neq 0$ and $\tilde{v}(0,\la)\neq 0$, then
$c_{1}=c_{2}
\frac{\tilde{v}(0,\la)}{\tilde{u}(0,\la)}$, and hence, $\D(E_{\lambda})=1$.\\
If $\tilde{u}(0,\lambda)=0$ and $\tilde{v}(0,\lambda)\neq 0$ (or
$\tilde{u}(0,\lambda)\neq 0$ and $\tilde{v}(0,\lambda)=0$), then
\begin{equation*} \phi(x,\la)= \left\{
\begin{array}{lll}
c_{1}\tilde{u}(x,\la),~&-1\leq x\leq 0,\\
0, ~&0\leq x\leq 1.\\
\end{array}
\right.
\end{equation*}
Since $\tilde{u}'(0,\la)\neq 0$ and $\tilde{v}'(0,\la)=0$, then the
last condition in Problem \eqref{an} is not satisfied, a
contradiction.\\
Now, if $\tilde{u}(0,\la)=\tilde{v}(0,\la)=0$, then from the last
condition in \eqref{an}, we have
\begin{equation}\label{ll}
c_{1}\sigma_{1}(0)\tilde{u}'(0)-c_{2}\sigma_{2}(0)\tilde{v}'(0)=0.
\end{equation}
Since $\tilde{u}'(0,\lambda) \neq 0$ and $\tilde{v}'(0,\lambda) \neq
0$, then
$c_{1}=c_{2}\frac{\sigma_{2}(0)\tilde{v}'(0,\lambda)}{\sigma_{1}(0)
\tilde{u}'(0,\lambda)}$,
whence $\D(E_{\lambda})=1$.\\
Since the operator $\mathcal{A}$ is self-adjoint, then the algebraic
multiplicity of $\la$ is equal to one.
\end{Proof}\\
Now, we enunciate our main result in this section:
\begin{theo}\label{Tg1}
There exists a subsequence of eigenvalues
$(\la_{\varphi(n)})_{n\geq1}$ which satisfies the asymptotic
estimates:\bean
&~&\sqrt{\la_{\varphi(n)+1}}-\sqrt{\la_{\varphi(n)}}=
\mathcal{O}(\dfrac{1}{n}),\label{g26}\\
%\la_{\varphi(n)+1}-\la_{\varphi(n)}=\mathcal{O}(\dfrac{1}{n}),\\
&~&\sqrt{\la_{n+2}}-\sqrt{\la_{n}}>2\delta, ~n\geq1,~\hbox{for
some}~\delta>0.\label{g27} \eean
\end{theo}
In order to prove this theorem, we establish some preliminary
results. Let\\
$\Gamma=\{\mu_{n}\}_{1}^{\infty}=\{\mu^{-}_{n}\}_{1}^{\infty}\bigcup\{\mu^{+}_{n}\}_{1}^{\infty}$,
where $\mu^{-}_{n}$ and $\mu^{+}_{n}$ are the eigenvalues of the
problems\begin{equation}\label{g1}  \left\{
\begin{array}{lll}
-(\sigma_{1}(x)y_{x}(x))_{x}+q_{1}(x)y(x)=\lambda\rho_{1}(x)y(x),~x\in (-1,\,0),\\
y(-1)=y(0)=0
\end{array}
\right.
\end{equation}and
 \begin{equation}\label{g2}
\left\{
\begin{array}{lll}
-(\sigma_{2}(x)y_{x}(x))_{x}+q_{2}(x)y(x)=\lambda\rho_{2}(x)y(x),~x\in (0,\,1),\\
y(0)=y(1)=0,
\end{array}
\right.
\end{equation}
respectively. It is known \cite[Chapter 1]{BI}, that the eigenvalues
$\mu^{-}_{j}$ and $\mu^{+}_{k}$ satisfy the asymptotic estimates
\bean\label{la} \left\{
 \begin{array}{ll}
   \sqrt{\mu^{-}_{j}}=\frac{j\pi}{\gamma_{1}}+\mathcal{O}(\frac{1}{j}), \\
   \sqrt{\mu^{+}_{k}}=\frac{k\pi}{\gamma_{2}}+\mathcal{O}(\frac{1}{k}),
\end{array}
 \right.
\eean where \bean
\gamma_{1}=\int_{-1}^{0}\sqrt{\frac{\rho_{1}(x)}{\sigma_{1}(x)}}dx~~and~~~
\gamma_{2}=\int_{0}^{1}\sqrt{\frac{\rho_{2}(x)}{\sigma_{2}(x)}}dx.
\label{gamm}\eean Obliviously, $\mu^{-}_{j}$ and $\mu^{+}_{k}$ can
coincide. Let
\begin{equation}\label{g34}
\Gamma^{*}=\{\mu_{n}\in\Gamma~\backslash~\mu^{-}_{j}=\mu^{+}_{k},~~for~some~j,~k\in\mathbb{N^{*}}\}.
\end{equation}
Note that if $\mu_{n}\in\Gamma^{*}$ (i.e.,
$\tilde{u}(0,\mu_n)=\tilde{v}(0,\mu_n)=0$), then $\mu_n$ is an
eigenvalue of both problems \eqref{g1} and \eqref{g2}. Clearly, all
the eigenvalues $\mu_{n}\in \Gamma$ can be arranged as follows:
\begin{equation*}
0<\mu_{1}\leq\mu_{2}\leq.......\leq\mu_{n}\leq.....\underset{n\rightarrow
\infty}{\longrightarrow}\infty.\end{equation*}
\begin{rem}\begin{itemize}
\item If the coefficients $\rho_i$, $\sigma_i$ and $q_i~(i=1,2)$ are symmetric,
then $\Gamma\equiv\Gamma^{*}$.
\item In the case of constant coefficients $\sigma_i\equiv1$, $q_i\equiv0
~(i=1,2)$, if
$\sqrt{\frac{\rho_1}{\rho_2}}\in\mathbb{R}\backslash\mathbb{Q}$,
then $\Gamma^{*}\equiv\emptyset.$
\end{itemize}
\end{rem}
We consider the following boundary value problem
\begin{equation}\label{g3}
\left\{
\begin{array}{lll}
-(\sigma_{1}(x)u_{x})_{x}+q_{1}(x)u=\lambda \rho_{1}(x)u,~~~x\in (-1,\,0),\\
-(\sigma_{2}(x)v_{x})_{x}+q_{2}(x)v=\lambda \rho_{2}(x)v,~~~x\in (0,\,1),\\
u(-1)=v(1)=0,\\
u(0)=v(0).\\
\end{array}
\right.
\end{equation}
It is clear that for $\la \in \mathbb{C}\backslash\Gamma$, the set
of solutions of Problem \eqref{g3} is a one-dimensional subspace
which is generated by a solution of the form
\begin{equation}\label{g4}
\widetilde{U}(x,\la)= \left\{
\begin{array}{lll}
\tilde{v}(0,\la)\tilde{u}(x,\la),&-1\leq x\leq 0,\\
\tilde{u}(0,\la)\tilde{v}(x,\la),&0\leq x\leq 1,\\
\end{array}
\right.
\end{equation}
where $\tilde{u}(x,\la)$ and $\tilde{v}(x,\la)$ are the solutions of
the initial value problems \eqref{w1} and \eqref{w2}, respectively.
Note that $\tilde{u}(0,\la)\neq0$ and $\tilde{v}(0,\la)\neq0$ for
$\la \in(\mu_{n},\mu_{n+1})$, since otherwise $\la$ would be an
eigenvalue of Problem \eqref{g1} or \eqref{g2}. Let us introduce the
variable complex function
$$F(\lambda)=\dfrac{\sigma_{1}(0)\widetilde{U}_{x}(0^-,\,\lambda)-
\sigma_{2}(0)\widetilde{U}_{x}(0^+,\,\lambda)}{\widetilde{U}(0,\,\lambda)},~\la\in\mathbb{C}\backslash\Gamma,$$
which can be rewritten in the form \be
F(\lambda)=\frac{\sigma_{1}(0)\tilde{v}(0)\tilde{u}_{x}(0^-)-
\sigma_{2}(0)\tilde{u}(0)\tilde{v}_{x}(0^+)}{\tilde{u}(0)\tilde{v}(0)},~\la\in\mathbb{C}\backslash\Gamma.\label{g5}\ee
It is known \cite[Chapter 1]{BI}, that $\tilde{u}(x,\la)$ and
$\tilde{v}(x,\la)$ are entire functions in $\la$ and continuous on
the intervals $[-1,0]$ and $[0,1]$, respectively. Therefore,
$F(\la)$ is a
 meromorphic function. We will show below that its zeros and poles coincide with the
eigenvalues of the regular problem \eqref{an} for $M=0$ and the
eigenvalues $\mu_{n}$, $n\geq1$, respectively. Moreover, the
solution of the equation \be F(\la)=M\la,\label{g14}\ee are the
eigenvalues $\la_{n}$, $n\geq 1$ of Problem \eqref{an}.
\begin{lemm}\label{Lg1}
The function $F(\lambda)$ is decreasing along the intervals
$(-\infty,\mu_{1})$ and\\$(\mu_{n}, \mu_{n+1})$, $n\geq1$ with
$\mu_{n}\neq\mu_{n+1}$. Furthermore, it decreases from $+\infty$ to
$-\infty$ in all intervals.
\end{lemm}
\begin{Proof}
Let $(\la, \la')\in (\mu_{n}, \mu_{n+1})\times(\mu_{n}, \mu_{n+1})$,
where $\la\neq\la'$ and $\widetilde{U}(x,\la)$,
$\widetilde{U}(x,\la')$ be the solutions of Problem \eqref{g3}.
Integrating by parts and taking into account the boundary conditions
in \eqref{g3}, yield \bean\label{g6} \left\{
 \begin{array}{ll}
    \displaystyle\int_{-1}^{0}\left( (q_{1}-\la \rho_{1})(x)\tilde{u}(x,\la)
    \tilde{u}(x,\la')+\sigma_{1}(x) \tilde{u}_{x}(x,\la)\tilde{u}_{x}(x,\la')\right)dx
     = \sigma_{1}(0)\tilde{u}(0,\la')\tilde{u}_{x}(0^-,\la), \\
    \displaystyle\int_{-1}^{0}\left( (q_{1}-\la' \rho_{1})(x)\tilde{u}
    (x,\la')\tilde{u}(x,\la)+\sigma_{1}(x) \tilde{u}_{x}(x,\la')\tilde{u}_{x}(x,\la)\right)dx
    =\sigma_{1}(0) \tilde{u}(0,\la)\tilde{u}_{x}(0^-,\la')
\end{array}
 \right.
\eean and \bean\label{g7} \left\{
 \begin{array}{ll}
    \displaystyle\int_{0}^{1}\left( (q_{2}-\la \rho_{2})(x)\tilde{v}(x,\la)\tilde{v}(x,\la')
    +\sigma_{2}(x) \tilde{v}_{x}(x,\la) \tilde{v}_{x}(x,\la')\right)dx
     = -\sigma_{2}(0)\tilde{v}(0,\la')\tilde{v}_{x}(0^+,\la), \\
    \displaystyle\int_{0}^{1}\left( (q_{2}-\la' \rho_{2})(x)\tilde{v}(x,\la')
    \tilde{v}(x,\la)+\sigma_{2}(x) \tilde{v}_{x}(x,\la') \tilde{v}_{x}(x,\la)\right)dx
    = -\sigma_{2}(0) \tilde{v}(0,\la)\tilde{v}_{x}(0^+,\la').
\end{array}
 \right.
\eean Subtracting the two equations of Systems \eqref{g6} and
\eqref{g7}, one gets\bean \left\{
 \begin{array}{ll}
    (\la'-\la)\displaystyle\int_{-1}^{0} \rho_{1}(x)\tilde{u}(x,\la)\tilde{u}(x,\la')dx
     = \sigma_{1}(0) (\tilde{u}(0,\la')\tilde{u}_{x}(0^-,\la)-\tilde{u}(0,\la)\tilde{u}_{x}(0^-,\la')), \\
     (\la'-\la)\displaystyle\int_{0}^{1} \rho_{2}(x)\tilde{v}(x,\la)\tilde{v}(x,\la')dx
     = \sigma_{2}(0)(\tilde{v}(0,\la)\tilde{v}_{x}(0^+,\la')-\tilde{v}(0,\la')\tilde{v}_{x}(0^+,\la)).
\end{array}
 \right.
\eean Hence \bean (\la-\la')\displaystyle\int_{-1}^{0}
\rho_{1}(x)\tilde{u}(x,\la)\tilde{u}(x,\la')dx
     &=& \sigma_{1}(0) (\tilde{u}(0,\la)(\tilde{u}_{x}(0^-,\la')-\tilde{u}_{x}(0^-,\la))\nonumber\\
     &&-\sigma_{1}(0)\tilde{u}_{x}(0^-,\la)(\tilde{u}(0,\la')-\tilde{u}(0,\la))~~~~~~
\eean and \bean (\la'-\la)\displaystyle\int_{0}^{1}
\rho_{2}(x)\tilde{v}(x,\la)\tilde{v}(x,\la')dx
&=& \sigma_{2}(0) (\tilde{v}(0,\la)(\tilde{v}_{x}(0^+,\la')-\tilde{v}_{x}(0^+,\la))\nonumber\\
     &&-\sigma_{2}(0)\tilde{v}_{x}(0^+,\la)(\tilde{v}(0,\la')-\tilde{v}(0,\la)).~~~~~~
\eean Passing to the limit as $\la' \rightarrow \la$, we get the
identities \bean\label{g16} \left\{
 \begin{array}{ll}
-\displaystyle\int_{-1}^{0} \rho_{1}(x)\tilde{u}^{2}(x,\la)dx
=\sigma_{1}(0)
 (\tilde{u}(0,\la)\dfrac{\partial \tilde{u}_{x}(0^-,\la)}{\partial\la}-\tilde{u}_{x}(0^-,\la)
 \dfrac{\partial \tilde{u}(0,\la)}{\partial\la}),\\
-\displaystyle\int_{0}^{1} \rho_{2}(x)\tilde{v}^{2}(x,\la)dx=
-\sigma_{2}(0) (\tilde{v}(0,\la)\dfrac{\partial
\tilde{v}_{x}(0^+,\la)}{\partial \la}-\tilde{v}_{x}(0^+,\la)\dfrac
{\partial \tilde{v}(0,\la)}{\partial \la}).
\end{array}
 \right.
\eean Dividing the first equation in \eqref{g16} by
$\tilde{u}^2(0,\la)$ and the second by $\tilde{v}^2(0,\la)$, it
follows \be\frac{\partial
F(\la)}{\partial\la}=-\dfrac{\tilde{v}^2(0,\la)\int_{-1}^{0}
\rho_{1}(x)\tilde{u}^{2}(x,\la)dx+ \tilde{u}^2(0,\la)\int_{0}^{1}
\rho_{2}(x)\tilde{v}^{2}(x,\la)dx}{\tilde{u}^2(0,\la)\tilde{v}^2(0,\la)}<0.\label{g17}\ee
In order to prove the second statement, we firstly establish the
asymptotic of $F(\la)$ as $\la\rightarrow-\infty$. It is known (e.g;
\cite[Chapter 2]{M} and \cite[Chapter 1]{BI}) for $\la\in\mathbb{C}$
and $|\la|\rightarrow \infty$ that \bean\label{g21}
\begin{cases}
\tilde{u}(x,\la)=a_{1}
(\rho_{1}(x)\sigma_{1}(x))^{-\frac{1}{4}}\dfrac{\sin(\sqrt{\la}\int_{-1}^{x}
\sqrt{\frac{\rho_{1}(t)}{\sigma_{1}(t)}})}{\sqrt{\la}}[1],\\
 \tilde{u}_{x}(x,\la)=
a_{1}(\rho_{1}(x))^{\frac{1}{4}}(\sigma_{1}(x))^{-\frac{3}{4}}\cos(\sqrt{\la}\int_{-1}^{x}
\sqrt{\frac{\rho_{1}(t)}{\sigma_{1}(t)}}dt)[1]
\end{cases}
\eean and \bean\label{g22}
\begin{cases}
\tilde{v}(x,\la)=a_{2}
(\rho_{2}(x)\sigma_{2}(x))^{-\frac{1}{4}}\dfrac{\sin(\sqrt{\la}\int_{x}^{1}
\sqrt{\frac{\rho_{2}(t)}{\sigma_{2}(t)}}dt)}{\sqrt{\la}}[1],\\
\tilde{v}_{x}(x,\la)=
-a_{2}(\rho_{2}(x))^{\frac{1}{4}}(\sigma_{2}(x))^{-\frac{3}{4}}\cos
(\sqrt{\la}\int_{x}^{1}
\sqrt{\frac{\rho_{2}(t)}{\sigma_{2}(t)}}dt)[1],
\end{cases}
\eean where $[1]=1+\mathcal{O}(\frac{1}{\sqrt{\la}})$, $\gamma_{i}$
$(i=1,2)$ are defined by \eqref{gamm} and the constants
\be\label{am27}
\begin{cases} a_1=\(\left(\rho_{1}(-1)\right)^{\frac{1}{4}}
\left(\sigma_{1}(-1)\right)^{-\frac{3}{4}}\)^{-1},\\
a_2=\(\left(\rho_{2}(1)\right)^{\frac{1}{4}}
\left(\sigma_{2}(1)\right)^{-\frac{3}{4}}\)^{-1},
\end{cases}\ee
are determined by the initial conditions in \eqref{w1} and
\eqref{w2}. By use of \eqref{g21} and \eqref{g22}, a straightforward
calculation gives the following asymptotic
\begin{equation*}
F(\la)\sim\sqrt{|\la|}\big(\sqrt{\rho_{1}(0)\sigma_{1}(0)}+\sqrt{\rho_{2}(0)\sigma_{2}(0)}\big),~as~\la\rightarrow-\infty.\end{equation*}
This implies that $\lim\limits_{\lambda\rightarrow
-\infty}F(\lambda)=+\infty$.\\
Now, we prove that \be
\lim_{\la\rightarrow\mu_{n}+0}F(\la)=+\infty,~~~~
\lim_{\la\rightarrow\mu_{n}-0}F(\la)=-\infty.\label{g18}\ee Suppose
that $\tilde{u}(0,\mu_{n})=\tilde{v}(0,\mu_{n})$=0. Let
$\la=\mu_{n}+\epsilon$, where $\epsilon$ is small enough. Therefore,
a simple calculation yields
 \bean F(\lambda)=&\dfrac{\sigma_{1}(0)\tilde{u}_{x}(0,\mu_{n})\frac{\partial \tilde{v}}{\partial\la}(0,\mu_{n})
-\sigma_{2}(0)\tilde{v}_{x}(0,\mu_{n})\frac{\partial
\tilde{u}}{\partial\la}(0,\mu_{n})}{\epsilon \frac{\partial
\tilde{u}}{\partial\la}(0,\mu_{n})\frac{\partial \tilde{v}}
{\partial\la}(0,\mu_{n})}+ \mathcal{O}(1),\nonumber \\
=&\frac{1}{\epsilon}\Big(\dfrac{\sigma_{1}(0)\tilde{u}_{x}(0,\mu_{n})}{
\frac{\partial \tilde{u}}{\partial\la}(0,\mu_{n})}-
\dfrac{\sigma_{2}(0) \tilde{v}_{x}(0,\mu_{n})}{\frac{\partial
\tilde{v}}{\partial\la}(0,\mu_{n})}\big)+o(1).\label{g13} \eean
Since $\mu_{n}$ is a simple eigenvalue of the two problems
\eqref{g1} and \eqref{g2}, then $\frac{\partial
\tilde{u}}{\partial\la}(0,\mu_{n}) \neq 0$ and
  $\frac{\partial \tilde{v}}{\partial\la}(0,\mu_{n})\neq0$.
Let us denote by
\begin{equation}\label{g1*}
F_{1}(\lambda)=\frac{\sigma_{1}(0)\tilde{u}_{x}(0,\la)}{\tilde{u}(0,\la)}
~~and~~
F_{2}(\lambda)=\frac{\sigma_{2}(0)\tilde{v}_{x}(0,\la)}{\tilde{v}(0,\la)}.
\end{equation}
It is easily seen that the eigenvalues $\mu^{-}_{j}$ and
$\mu^{+}_{k}$ $(j\geq1,~k\geq1)$ are the poles of $F_{1}(\lambda)$
and $F_{2}(\lambda)$, respectively. In view of the proof of
Proposition 4 in \cite{JA}, by use of the Mittag-Lefleur theorem
\cite[Chapter 4]{Hur}, $F_{1}(\la)$ has the following decomposition
$$F_{1}(\la)= \displaystyle\sum_{j\geq1}
(\dfrac{\la}{\mu^{-}_{j}})\dfrac{c_{j}}{\la-\mu^{-}_{j}},$$ where
$c_{j}$ are the residuals of $F_{1}(\la)$ at the poles
$\mu^{-}_{j}$, $j\geq1$. It is known that the residuals of
$F_{1}(\lambda)$ and $F_{2}(\lambda)$ (at the poles $\mu^{-}_{j}$
and $\mu^{+}_{k}$, respectively) are given by
$$c_{j}=\frac{\sigma_{1}(0)\tilde{u}_{x}(0,\mu^{-}_{j})}
{\frac{\partial \tilde{u}
}{\partial\la}(0,\mu^{-}_{j})},~~c_{k}^{'}= \frac{\sigma_{2}
(0)\tilde{v}_{x}(0,\mu^{+}_{k})}{\frac{\partial
\tilde{v}}{\partial\la}(0,\mu^{+}_{k})}.$$ According again to the
proof of Proposition 4 in \cite{JA}, we have $c_{j}>0$, $j\geq1$. By
the change of variables $s=-x$, if $\tilde{v}(x,\la)$ is a solution
of Problem \eqref{w2}, then $\tilde{v}(-s,\la)$ is a solution of
Problem \eqref{w1} (where the coefficients are replaced by
$\rho_{2}(-s)$, $\sigma_{2}(-s)$ and $q_{2}(-s)$), and hence,
$F_{2}(\la)=-\frac{\sigma_{2}(0)\tilde{v}_{s}(0,\la)}{\tilde{v}(0,\la)}$.
Therefore, $c'_{k}<0$ for all $k\geq1$ and
$$\dfrac{\sigma_{1}(0)\tilde{u}_{x}(0,\mu_{n})}{\frac{\partial
\tilde{u}}{\partial\la}(0,\mu_{n})}- \dfrac{\sigma_{2}(0)
\tilde{v}_{x}(0,\mu_{n})}{ \frac{\partial
\tilde{v}}{\partial\la}(0,\mu_{n})}>0,~n\geq1.$$Passing to the limit
as $\epsilon\rightarrow 0$ in \eqref{g13}, we obtain the
first limit in \eqref{g18}. Analogously, we can prove the second limit of \eqref{g18}.\\
Now, if $\tilde{u}(0,\mu_{n})=0$ and $\tilde{v}(0,\mu_{n})\neq0$ (or
$\tilde{v}(0,\mu_{n})=0$ and $\tilde{u}(0,\mu_{n})\neq0$), then
$\sigma_{1}(0)\tilde{u}_{x}(0,\mu_{n})\neq0$ (or
$\sigma_{2}(0)\tilde{v}_{x}(0,\mu_{n})\neq0$), and hence, we arrive
to the same conclusion.
\end{Proof}\\
From the last lemma, it is clear that the poles of $F(\lambda)$
coincide with the eigenvalues $\mu_{n}$, $n\geq 1$. Furthermore, it
follows the following interpolation formulas between the eigenvalues
$\la_{n}$, $\mu_{n}$ and those of the regular problem \eqref{an} for
$M=0$.
\begin{prop}\label{cg1}
 Let $\la'_{n}$, $n\geq1$ denote the eigenvalues of the regular
 problem \eqref{an} for
 $M=0$. If $\mu_{n} \neq \mu_{n+1}$, then
 \begin{equation}\label{lam}
\la_{1}<\la'_{1}<\mu_{1}~and~
\mu_{n}<\la_{n+1}<\la'_{n+1}<\mu_{n+1},~n\geq1
\end{equation}
and if $\mu_{n}=\mu_{n+1}$, then
$\mu_{n}=\lambda_{n+1}=\lambda'_{n+1}$. Furthermore the eigenvalues
$(\lambda_{n})_{n\geq1}$ satisfy the asymptotic
\begin{equation}\label{lbd}
\lambda_{n} \sim
\left(\frac{n\pi}{\gamma_{1}+\gamma_{2}}\right)^{2},
\end{equation}
where $\gamma_{1}$ and $\gamma_{2}$ are defined by \eqref{gamm}.
\end{prop}
\begin{Proof} According to Lemma \ref{Lg1}, $F(\la)$ is a decreasing function from $+\infty$ to $-\infty$ along each of intervals
 $(-\infty, \mu_{1})$ and $(\mu_{n}, \mu_{n+1})$, $n\geq1$. Hence, Equation \eqref{g14} has exactly
 one zero in each of these intervals. Moreover, the equation $F(\la)=0$  has exactly one
 zero in each of these intervals. It is clear that these zeros (denoted by $\la'_{n}$)
 are the eigenvalues of the regular problem \eqref{an} for $M=0$. Consequently, the
 interpolation formulas \eqref{lam}
 are simple deductions from the graphs of the functions $F(\la)$ and
 $M\la$.\\
It is known (e.g., \cite{FVAA} and \cite[Chapter 6.7]{IGK}), that
the eigenvalues $\la'_{n}$ satisfy the asymptotes
\be\label{g23}{\la'_{n}}=\(\frac{n\pi}{\int_{-1}^{1}
\sqrt{\frac{\rho(x)}{\sigma(x)}}dx}\)^2+\mathcal{O}(1),\ee where
\bea \rho(x)=\left\{
 \begin{array}{ll}
   \rho_{1}(x),~~x\in[-1,0], \\
   \rho_{2}(x),~~x\in[0,1],
\end{array}
 \right.
~~\mbox{and}~~ \sigma(x)=\left\{
 \begin{array}{ll}
   \sigma_{1}(x),~~x\in[-1,0], \\
   \sigma_{2}(x),~~x\in[0,1],
\end{array}
 \right.
. \eea Since ${\la'_{n-1}}<{\la_{n}}\leq{\la'_{n}}$, then
\eqref{lbd} follows from \eqref{g23}.
\end{Proof}\\
\begin{lemm}\label{Lem1} Let $(\mu_{n_{k}})$ be a subsequence of
$(\mu_{n})$ such that
\begin{equation}\label{mk}
\left|\sqrt{\mu_{n_{k}}}-\sqrt{\mu_{n_{k}-1}}\right|\rightarrow 0,
~\hbox{as}~k\rightarrow \infty.
\end{equation}
\begin{description}
  \item[(i$)$]
If $\mu_{n_{k}}$ is an eigenvalue of Problem \eqref{g1} (or of
Problem \eqref{g2}), then $\mu_{n_{k}-1}$ is an eigenvalue of
Problem \eqref{g2} (or of Problem \eqref{g1}).\\Furthermore,
$\sqrt{\mu_{n_{k}+1}}-\sqrt{\mu_{n_{k}}}>\delta_{1}$ and
$\sqrt{\mu_{n_{k}-1}}-\sqrt{\mu_{n_{k}-2}}>\delta_{2}$ for some
$\delta_{1}$, $\delta_{2}>0$.
\item[(ii$)$] If for some $C>0$,
\begin{equation}\label{am3}
\sqrt{\la_{n_{k}}}\left|\tilde{u}(0,\la_{n_{k}})\right| \leq
\frac{C}{k} ~(\hbox{or}~
\sqrt{\la_{n_{k}}}\left|\tilde{v}(0,\la_{n_{k}})\right| \leq
\frac{C}{k}), ~\hbox{for large}~k\in\mathbb{N}^*,
\end{equation}
\end{description}
then
\begin{equation}\label{am15}
\left|\tilde{u}(0,\la_{n_{k}})\right|\sim
\left|\tilde{v}(0,\la_{n_{k}})\right|,~\hbox{as}~k\rightarrow
\infty.
\end{equation}
\end{lemm}
\begin{Proof} \begin{description}
  \item[(i$)$] It is clear from the asymptotes \eqref{la} and
\eqref{mk}, that the eigenvalues $\mu_{n_{k}}$ and $\mu_{n_{k}-1}$
are of different types i.e., if $\mu_{n_{k}}$ is an eigenvalue of
Problem \eqref{g1}, then $\mu_{n_{k}-1}$ is an eigenvalue of Problem
\eqref{g2}, or conversely.\\Moreover,
$\sqrt{\mu_{n_{k}+1}}-\sqrt{\mu_{n_{k}}}>\delta_{1}$ and
$\sqrt{\mu_{n_{k}-1}}-\sqrt{\mu_{n_{k}-2}}>\delta_{2}$ for some
$\delta_{1},~\delta_{2}>0$.
%suppose that $\sqrt{\mu_{n_{k}+1}}-\sqrt{\mu_{n_{k}}}\rightarrow 0$
%as $k\rightarrow \infty$. Taking into account that $\mu_{n_{k}}$ is
%an eigenvalue of Problem \eqref{g1}, as above, we have
%$\mu_{n_{k}+1}$ is an eigenvalue of Problem \eqref{g2}. These imply
%that $\sqrt{\mu_{n_{k}+1}}-\sqrt{\mu_{n_{k}-1}}\rightarrow 0$ as
%$k\rightarrow \infty$, with $\mu_{n_{k}+1}$ and $\mu_{n_{k}-1}$ are
%eigenvalues of Problem \eqref{g2}. This is in contradiction with
%\eqref{la}. Similarly, we prove that
%$\sqrt{\mu_{n_{k}-1}}-\sqrt{\mu_{n_{k}-2}}>\delta_{2}$ for some
%$\delta_{2}>0$.
\item[(ii$)$] We may first assume for large
$k\in\mathbb{N}^*$, that
$\sqrt{\la_{n_{k}}}\left|\tilde{u}(0,\la_{n_{k}})\right| <
\frac{C}{k}, C>0$. Using the expression \eqref{g5} of $F\(\la\)$ and
the characteristic equation \eqref{g14},
\begin{equation}\label{eq}
\left|\frac{\tilde{v}(0,\la_{n_{k}})}{\tilde{u}(0,\la_{n_{k}})}\right|
=\left|\frac{\sigma_{2}(0)
\tilde{v}_{x}(0,\la_{n_{k}})}{\sigma_{1}(0)
\tilde{u}_{x}(0,\la_{n_{k}})-M\la_{n_{k}}\tilde{u}(0,\la_{n_{k}})}\right|,
\end{equation}
and hence, by \eqref{g21} and \eqref{g22}, one has
\begin{equation}\label{bb3}
\left|\frac{\tilde{v}(0,\la_{n_{k}})}{\tilde{u}(0,\la_{n_{k}})}\right|=
\left|\frac{a_{2}(\rho_{2}(0)\sigma_{2}(0))^{\frac{1}{4}}
\cos(\sqrt{\la_{n_{k}}}\gamma_{2})[1]
}{a_{1}\left((\rho_{1}(0)\sigma_{1}(0))^{\frac{1}{4}}
\cos(\sqrt{\la_{n_{k}}}\gamma_{1})- M
\sqrt{\la_{n_{k}}}(\rho_{1}(0)\sigma_{1}(0))^{\frac{-1}{4}}
\sin(\sqrt{\la_{n_{k}}}\gamma_1)\right)[1]} \right|,
\end{equation} the constants $a_{i}~(i=1,2)$ are defined by \eqref{am27} and
$[1]=1+\frac{1}{\sqrt{\la_{n_{k}}}}$. Under the above assumption
together with \eqref{lbd},
\begin{equation}\label{lam1}
\sqrt{\la_{n_{k}}}\left|\sin\(\sqrt{\la_{n_{k}}}\gamma_1\)
\right|[1]\rightarrow0
\hbox{ as } k\rightarrow\infty.
\end{equation}
In view of assertion $(i)$, we may assume that $\mu_{n_{k}-1}$ is an
eigenvalue of Problem \eqref{g1} and $\mu_{n_{k}}$ is an eigenvalue
of Problem \eqref{g2} (the other case can be handled in a same way).
We set
\begin{equation}\label{delta}
\delta_{n_{k}}^{-}=\sqrt{\la_{n_{k}}}-\sqrt{\mu_{n_{k}-1}}~\hbox{
and} ~\delta_{n_{k}}^{+}=\sqrt{\mu_{n_{k}}}-\sqrt{\la_{n_{k}}}.
\end{equation}
Hence, by \eqref{lam} and \eqref{mk},
\begin{equation}\label{delta1}
\delta_{n_{k}}^{+}\rightarrow 0~\hbox{and}~
\delta_{n_{k}}^{-}\rightarrow 0,~\hbox{as}~k\rightarrow \infty.
\end{equation}
In view of \eqref{la}, \eqref{delta} and \eqref{delta1}, we obtain
\begin{equation*}
\left|\sin(\sqrt{\la_{n_{k}}}\gamma_2)\right|\rightarrow 0,
~\hbox{as}~ k\rightarrow\infty.
\end{equation*}
From this and \eqref{lam1}, we have
$\cos(\sqrt{\la_{n_{k}}}\gamma_1)=\mathcal{O}(1)$ and
$\cos(\sqrt{\la_{n_{k}}}\gamma_2)=\mathcal{O}(1)$. Therefore, by
\eqref{bb3},
\bea\left|\frac{\tilde{v}\(0,\la_{n_{k}}\)}{\tilde{u}\(0,\la_{n_{k}}\)}\right|\sim
\frac{a_{2}}{a_{1}}\left(\frac{\rho_{2}(0)\sigma_{2}(0)}
{\rho_{1}(0)\sigma_{1}(0)}\right)^{\frac{1}{4}}, ~~\hbox{for $k$
large enough}.\eea Now, we assume for large $k\in\mathbb{N}^*$, that
\begin{equation}\label{ame}
\sqrt{\la_{n_{k}}}\left|\tilde{u}(0,\la_{n_{k}})\right|=
\frac{C}{k}+o(\frac{1}{k}), ~C>0.
\end{equation}
In the expression \eqref{eq}, suppose that
\begin{equation}\label{bb}
\sigma_{1}(0)\tilde{u}_{x}(0,\la_{n_{k}})-M
\la_{n_{k}}\tilde{u}(0,\la_{n_{k}})\rightarrow0, \hbox{ as }
k\rightarrow\infty.
\end{equation}
Substituting \eqref{g21} and \eqref{ame} into \eqref{bb}, one
obtains
\begin{equation*}
M\sqrt{\la_{n_{k}}}\frac{C}{k}=
(\rho_{1}(0)\sigma_{1}(0))^{\frac{1}{2}}+\mathcal{O}(\frac{1}{k}).
\end{equation*}
%In view of \eqref{lam}, we put
%\begin{equation}\label{delta}
%\delta_{n}^{-}=\sqrt{\la_{n_{k}}}-\sqrt{\mu_{n_{k}-1}}~\hbox{ and}
%~\delta_{n}^{+}=\sqrt{\mu_{n}}-\sqrt{\la_{n_{k}}}.
%\end{equation}
%Hence, by \eqref{mk}, $\delta_{n}^{+}\rightarrow 0$ and
%$\delta_{n}^{-}\rightarrow 0$.
Using this and \eqref{la} together with \eqref{delta} and
\eqref{delta1}, one gets
\begin{equation*}
M\frac{C\pi}{\gamma_{1}}= (\rho_{1}(0)\sigma_{1}(0))^{\frac{1}{2}}
~\hbox{ and}~ M\frac{C\pi}{\gamma_{2}}=
(\rho_{1}(0)\sigma_{1}(0))^{\frac{1}{2}}.
\end{equation*}
Then if $\gamma_{1}\neq \gamma_{2}$, this is in contradiction with
\eqref{bb}. If $\gamma_{1}=\gamma_{2}$, it is easily seen from
\eqref{g21} and \eqref{g22}, that \eqref{am15} holds. This completes
the proof.
\end{description}
\end{Proof}\\
%From
%\eqref{g21}, \eqref{g22} and \eqref{am3}, one has \begin{equation*}
%\left|\sin(\sqrt{\la_{n_{k}}}\gamma_{1})\right|+
%\left|\sin(\sqrt{\la_{n_{k}}}\gamma_{2})\right|\rightarrow 0, \hbox{
%as } k\rightarrow\infty.
%\end{equation*}
%From this and \eqref{la}, it follows that $\delta_{n}^{+}\rightarrow
%0$ and $\delta_{n}^{-}\rightarrow 0$.
%As we shall see in the proof of Theorem \ref{Tg1} and the
%controllability section, the following proposition is crucial.
%\begin{prop} \label{prop4} For some $\tau>0$, let \be{\Omega}
%=\left\{n\in \mathbb{N}^{*}~:~ {{\mu_{n}}}-{{\mu _{n-1}}}\geq\tau
%\right\}. \label{am4}\ee Then
% \be
%\sqrt{{\la_{n}}}-\sqrt{{\mu _{n-1}}}= \mathcal{O}(\frac{1}{n})
%~~\hbox{ for }~~ n\in \Omega. \label{g19}\ee
%\end{prop}
%\begin{rem}\label{Rem2}  From the asymptotics \eqref{la},
%it is clear that for large $n\in\N^*$ if
%$\left|\sqrt{\mu_{n}}-\sqrt{\mu _{n-1}}\right|$ close to zero, then
%both $n-1, n+1\in\Omega$. In this case, $\mu_{n-1}$ is an eigenvalue
%of Problem \eqref{g1} and $\mu_{n}$ is an eigenvalue of Problem
%\eqref{g2} (or conversely).
%\end{rem}
%\begin{Proof}
% This completes the proof.
%\end{Proof}\\
We are now ready to prove Theorem \ref{Tg1}.\\
\begin{Proof} First, we establish the asymptotic estimate
 \be
\sqrt{{\la_{n}}}-\sqrt{{\mu _{n-1}}}=
\mathcal{O}(\frac{1}{n}),~\hbox{ for }~ n\in \Omega, \label{g19}\ee
where for some $\tau>0$, \be{\Omega} =\left\{n\in \mathbb{N}^{*}~:~
{{\mu_{n}}}-{{\mu _{n-1}}}\geq\tau \right\}.\label{am4}\ee
%From the
%asymptotics \eqref{la}, it is clear that for large $n\in\N^*$ if
%$\left|\sqrt{\mu_{n}}-\sqrt{\mu _{n-1}}\right|$ close to zero, then
%both $n-1, n+1\in\Omega$. In this case, $\mu_{n-1}$ is an eigenvalue
%Problem \eqref{g1} and $\mu_{n}$ is an eigenvalue of
%Problem \eqref{g2} (or conversely).\\
First of all, let us recall from the asymptotics \eqref{la} that for
large $n\in\Omega$, we have \be \sqrt{\mu_{n}}-\sqrt{\mu
_{n-1}}\geq\frac{\tau_0}{n}, \label{am24}\ee for some $\tau_0>0$.
For $\xi \in (\sqrt{\mu_{n-1}}, \sqrt{\la_{n}}]$ we put
$G(\xi)=\dfrac{1}{F(\xi^2)}$, where $F$ is defined by \eqref{g5}. In
view of Proposition \ref{cg1}, $(\sqrt{\mu_{n-1}},
\sqrt{\la_{n}}]\subset(\sqrt{\mu_{n-1}},\sqrt{\la'_{n}})$. Hence,
$\sigma_{1}(0)\tilde{v}(0,\xi^2)\tilde{u}_{x}(0,\xi^2)-\sigma_{2}(0)\tilde{
u}(0,\xi^2)\tilde{v}_{x}(0,\xi^2)\neq0$ for $\la \in
(\sqrt{\mu_{n-1}}, \sqrt{\la_{n}}]$ and this implies that $G(\xi)$
is well-defined in this interval. By use of the mean value theorem
on the interval
 $[\sqrt{\mu_{n-1}^{\epsilon}}, \sqrt{\la_{n}}]$,
 (where $\sqrt{\mu_{n-1}^{\epsilon}}=\sqrt{\mu_{n-1}}+\epsilon$ for enough small
 $\epsilon>0$), we have
$$G(\sqrt{\la_{n}})-G(\sqrt{\mu_{n-1}^{\epsilon}})=
(\sqrt{\la_{n}}-\sqrt{\mu_{n-1}^{\epsilon}})\frac{\partial
G(\xi)}{\partial \xi}\mid_{\xi=\alpha_{n}},$$ for some $
\alpha_{n}\in (\sqrt{\mu_{n-1}^{\epsilon}}, \sqrt{\la_{n}})$. From
\eqref{g18}, one has
$\sqrt{\mu_{n-1}^{\epsilon}}\rightarrow\sqrt{\mu_{n-1}}$ and\\
$G(\sqrt{\mu_{n-1}^{\epsilon}})\rightarrow 0 $ as
$\epsilon\rightarrow0$. Then by using Equation \eqref{g14} and the
expression \eqref{g17} of $\frac{\partial F(\la)}{\partial \la}$, we
obtain \begin{equation}\label{g20}
\left(\sqrt{{\la}_{n}}-\sqrt{{\mu}_{n-1}}\right)=\left(\frac{\partial
G}{\partial\xi}(\alpha_{n})\right)^{-1}\left(\dfrac{1}{M\la_{n}}\right),\end{equation}
where \begin{equation*} \left(\frac{\partial
G}{\partial\xi}(\alpha_{n})\right)^{-1}=\dfrac{\Big(\sigma_{1}(0)\tilde
{v}(0,\alpha_{n}^2)\tilde{u}_{x}(0,\alpha_{n}^2)-
\sigma_{2}(0)\tilde{u}(0,\alpha_{n}^2)\tilde{v}_{x}(0,\alpha_{n}^2)
\Big)^{2}}{2\alpha_{n}\(\tilde{v}^2(0,\alpha_{n}^2)\int_{-1}^{0}
\rho_{1}(s)\tilde{u}^{2}(s,\alpha_{n}^2)ds+
\tilde{u}^2(0,\alpha_{n}^2)\int_{0}^{1}
\rho_{2}(s)\tilde{v}^{2}(s,\alpha_{n}^2)ds\)}.\end{equation*} Taking
into account the asymptotes \eqref{g21} and \eqref{g22}, one has
\bean \alpha_{n}
\Big(\sigma_{1}(0)\tilde{v}(0,\alpha_{n}^2)\tilde{u}_{x}(0,\alpha_{n}^2)-
\sigma_{2}(0)\tilde{u}(0,\alpha_{n}^2)\tilde{v}_{x}(0,\alpha_{n}^2)\Big)
=a_{1}a_{2}\Psi\(\alpha_{n}\), \label{cvk}\eean where
\begin{equation}\label{am5}
\Psi\(\xi\)=\left(\left(\frac{\rho_1(0)\sigma_1(0)}{\rho_2(0)\sigma_2(0)}\right)^{\frac{1}{4}}
\sin(\xi\gamma_2)\cos
(\xi\gamma_1)+\(\frac{\rho_1(0)\sigma_1(0)}{\rho_2(0)\sigma_2(0)}\)
^{-\frac{1}{4}}\sin(\xi\gamma_1)\cos (\xi\gamma_2)\right)[1],
\end{equation} $[1]=1+\mathcal{O}(\frac{1}{\xi})$. Hence,
by the Riemann-Lebesgue Lemma, \begin{equation}\label{cvs}
\frac{1}{\alpha_{n}} \(\frac{\partial
G}{\partial\xi}(\alpha_{n})\)^{-1}=
\frac{2\(\Psi\(\alpha_{n}\)\)^2}{\Big({\gamma_1}
{\(\rho_1(0)\sigma_1(0)\)^{-\frac{1}{2}}}\sin^2
\left({\alpha_{n}}\gamma_1\)+{\gamma_2}{\(\rho_2(0)
\sigma_2(0)\)^{-\frac{1}{2}}}\sin^2
\left({\alpha_{n}}\gamma_2\)\Big)[1]}. \end{equation} It is clear
that if $\left|\Psi\(\alpha_{n}\)\right|>\tau_1$ for some $\tau_1>0$
and all $n\geq 1$, then\be \frac{1}{\alpha_{n}}\(\frac{\partial
G}{\partial\xi}(\alpha_{n})\)^{-1}\asymp 1 ,\label{cva2}\ee and
hence, by \eqref{g20}, we get
\bean\(\sqrt{{{\la}_{n}}}-\sqrt{{{\mu}_{n-1}}}\)\asymp
\dfrac{\alpha_{n}}{\la_{n}}.\label{g20*}\eean Since
$\sqrt{{\la_{n-1}}}<{\alpha_{n}}<\sqrt{{\la_{n}}}$, then from
\eqref{lbd}, one has \bean\alpha_{n}\sim
\dfrac{n\pi}{\gamma_{1}+\gamma_{2}},~~as~~n\rightarrow\infty.\label{g24}
\eean Therefore, \eqref{g19} follows from the asymptotics
\eqref{lbd}, \eqref{g20*} and \eqref{g24}.\\ Now, we assume that
$\left|\Psi\(\alpha_{n_k}\)\right|\rightarrow0 \hbox{ as }
k\rightarrow\infty$, for some
 subsequence $(\alpha_{n_k})$ of $(\alpha_{n})$. Then from \eqref{cvk},
 $\alpha_{n_k}$ tends to the square root of an eigenvalue
of the regular problem \eqref{an} for $M=0$.\\In view of Proposition
\ref{cg1},
$\sqrt{{\la_{n_k-1}'}}<\sqrt{\mu_{n_k-1}}<{\alpha_{n_k}}<\sqrt{{\la_{n_k}}}
<\sqrt{{\la_{n_k}'}}$, and hence, $\alpha_{n_k}$ tends to
$\sqrt{\la'_{n_k-1}}$ or $\sqrt{\la'_{n_k}}$. According to
\eqref{am24},
two cases must be examined.\\
{\bf Case $1-$} Assume that for large $k\geq 1$ and $n_{k} \in
\Omega$,
 $\left|\sqrt{\mu_{n_k}}-\sqrt{\mu
_{n_k-1}}\right|\geq{\tilde{\tau}_{0}},$ for some
$\tilde{\tau}_{0}>0$. Under this assumption, we prove \eqref{cva2}.
Firstly, if $\left|{\alpha_{n_k}}-\sqrt{{\la_{n_k}'}}\right|$ tends
to zero, then \be
\left|\sqrt{\la_{n_k}}-\sqrt{{\la_{n_k}'}}\right|\rightarrow0,~\hbox{i.e.},~
\Psi\(\sqrt{\la_{n_k}}\)\rightarrow0 \hbox{ as }
k\rightarrow\infty.\label{g2*}\ee On the other hand by \eqref{g21},
\eqref{g22} and the characteristic equation \eqref{g14}, one gets
$$\Psi\(\sqrt{\la_{n_k}}\)= M
\(\rho_1(0)\sigma_1(0)\)^{-\frac{1}{4}}\(\rho_2(0)\sigma_2(0)\)^{-\frac{1}{4}}
\sqrt{\la_{n_k}}
\(\sin(\sqrt{\la_{n_k}}\gamma_1)\sin(\sqrt{\la_{n_k}}\gamma_2)\)[1],$$
 and hence, by \eqref{g2*},
$\sin(\sqrt{\la_{n_k}}\gamma_1)\[1\]\rightarrow0$ or
$\sin(\sqrt{\la_{n_k}}\gamma_2)\[1\]\rightarrow0$, as
$k\rightarrow\infty$ . Again by \eqref{g2*} together with
\eqref{am5}, we have
\begin{equation*}
\(\left|\sin(\sqrt{\la_{n_k}}\gamma_1)\right|+
\left|\sin(\sqrt{\la_{n_k}}\gamma_2)\right|\)[1]\rightarrow0,~\hbox{as}~
k\rightarrow\infty.\end{equation*} This means that
$\sqrt{\la_{n_k}}$ tends simultaneously to the square root of
eigenvalues of the two Problems \eqref{g1} and \eqref{g2}. Since
$\left|\sqrt{\mu_{n_k}}-\sqrt{\mu
_{n_k-1}}\right|\geq{\tilde{\tau}_{0}}$, then in view of assertion
$(i)$ of Lemma \ref{Lem1}, we have
 $\left|\sqrt{{\la_{n_k}}}-\sqrt{\mu_{n_k-1}}\right|$ tends to zero and
${\sqrt{\mu_{n_k-2}}}$ close to $\sqrt{\mu_{n_k-1}}$ (or
$\left|\sqrt{{\la_{n_k}}}-\sqrt{\mu_{n_k}}\right|\rightarrow 0$ and
$\left|\sqrt{{\mu_{n_{k}+1}}}-\sqrt{\mu_{n_k}}\right|\rightarrow
0$). According to Proposition \ref{cg1},
$\sqrt{\la'_{n_k-1}}\in\(\sqrt{{\mu_{n_k-2}}},\sqrt{\mu_{n_k-1}}\)$,
thus by \eqref{g2*}, one has
$$\left|\sqrt{\la'_{n_k}}-\sqrt{\la'_{n_k-1}}\right|\rightarrow0 ,~~as~~
k\rightarrow\infty.$$ This is in contradiction with the asymptotic
estimate \eqref{g23}.\\
Now, if $\left|{\alpha_{n_k}}-\sqrt{{\la_{n_k-1}'}}\right|$ converge
to zero, then $\left|{\alpha_{n_k}}-\sqrt{{\mu_{n_k-1}}}\right|$
tends to zero, and hence,
 $\sin({\alpha_{n_k}}\ga_1)\[1\]\rightarrow0$ or
$\sin({\alpha_{n_k}}\ga_2)\[1\]\rightarrow0$ as $
k\rightarrow\infty$. As above, since\\
$\Psi\(\alpha_{n_k}\)\rightarrow0~\hbox{as}~k\rightarrow\infty$,
then by \eqref{am5},
\begin{equation}\label{am}
\(|\sin({\alpha_{n_k}}\ga_1)|+
|\sin({\alpha_{n_k}}\ga_2)|\)\[1\]\rightarrow0,~\hbox{as}~
k\rightarrow\infty.\end{equation} As before in view of the
assumption
$\left|\sqrt{{\mu_{n_k}}}-\sqrt{\mu_{n_k-1}}\right|>\tilde{\tau}_0$,
\eqref{am} implies \be
\left|{\alpha_{n_k}}-\sqrt{{\mu_{n_k-1}}}\right|\rightarrow0 \hbox{
and }
\left|\sqrt{{\mu_{n_k-1}}}-\sqrt{\mu_{n_k-2}}\right|\rightarrow0,
\hbox{ as } k\rightarrow\infty.\label{am6}\ee Let us recall from the
expression
 \eqref{g5} of $F(\la)$ that $$F(\xi^2)= F_1(\xi^2)-F_2(\xi^2),$$
 where $F_1$ and $F_2$ are defined by \eqref{g1*}.
 This means that $\(\frac{\partial
 G}{\partial\xi}(\xi)\)^{-1}$ can be expressed in the form
\begin{equation}\label{gap1+}
\(\frac{\partial
G}{\partial\xi}(\xi)\)^{-1}=\frac{-\(F_1(\xi^2)-F_2(\xi^2)\)^2}
{2\xi\(\frac{\partial F_1}{\partial\xi}(\xi^2)-\frac{\partial
F_2}{\partial\xi}(\xi^2)\)}, \hbox{ for }\xi \in
\(\sqrt{\mu_{n_k-1}}, \sqrt{\la_{n_k}}\].\end{equation} Similarly to
the proof of Lemma \ref{Lg1}, we can prove that $F_1\(\xi^2\)$ is a
decreasing function from $+\infty$ to $-\infty$ in all the intervals
 $(-\infty, \sqrt{\mu_{n}^{-}})$ and $(\sqrt{\mu_{n}^{-}}, \sqrt{\mu_{n+1}^{-}})$,
 while, $F_2\(\xi^2\)$ is an increasing function from $-\infty$
to $+\infty$ along each of the intervals
 $(-\infty, \sqrt{\mu_{n}^{+}})$ and $(\sqrt{\mu_{n}^{+}},
 \sqrt{\mu_{n+1}^{+}})$, $n\geq1$.
By \eqref{am6}, $\sqrt{\mu_{n_k-2}}$ close to $\sqrt{\mu_{n_k-1}}$,
then in view of assertion $(i)$ of Lemma \ref{Lem1}, we may assume
that $\mu_{n_k-2}$ and $\mu_{n_k-1}$ are eigenvalues of Problems
\eqref{g1} and \eqref{g2}, respectively (the other case can be
handled in a same way), i.e.,\be
\lim_{\xi\rightarrow\sqrt{\mu_{n_k-2}}+0}F_1(\xi^2)=+\infty,~~~~
\lim_{\xi\rightarrow\sqrt{\mu_{n_k-1}}+0}F_2(\xi^2)=-\infty.\label{g18+}\ee
Therefore, by \eqref{am6} and \eqref{g18+}, for each $\alpha_{n_k}$
 in a sufficiently small right neighborhood of
$\sqrt{\mu_{n_k-1}}$, one gets \be
F_1(\alpha_{n_k}^2)F_2(\alpha_{n_k}^2)<0.\label{am25}\ee It is clear
that the above together with \eqref{gap1+} and \eqref{am25} imply
that \be\(\frac{\partial G}{\partial\xi}(\alpha_{n_k})\)^{-1}\asymp
\frac{\(F_1(\alpha_{n_k}^2)\)^2+\(F_2(\alpha_{n_k}^2)\)^2}
{\alpha_{n_k}\(\frac{\partial
F_2}{\partial\xi}(\alpha_{n_k}^2)-\frac{\partial
F_1}{\partial\xi}(\alpha_{n_k}^2)\)}.\label{am28}\ee Thus, combining
\eqref{cvs} with \eqref{am28}, we obtain $$ \frac{1}{\alpha_{n_k}}
\(\frac{\partial G}{\partial\xi}(\alpha_{n_k})\)^{-1}\asymp
\dfrac{\(\sin^2({\alpha_{n_k}}\ga_2)\cos^2
({\alpha_{n_k}}\ga_1)+\sin^2({\alpha_{n_k}}\ga_1)\cos^2
({\alpha_{n_k}}\ga_2)\)\[1\]}{\Big(\sin^2
\left({\alpha_{n_k}}\ga_1\)+\sin^2\left({\alpha_{n_k}}
\ga_2\right)\Big)\[1\]}%\no
$$ and hence, by \eqref{am}, the desired
 estimate \eqref{cva2} follows.\\
{\bf Case $2$-} Now, suppose that  $\sqrt{\mu_{n_k-1}}$ close to
 $\sqrt{\mu_{n_k}}$. This implies that for\\
 $\alpha_{n_k}\in\(\sqrt{\mu_{n_k-1}},~\sqrt{\mu_{n_k}}\)$,
 \begin{equation*}\(\left|\sin({\alpha_{n_k}}\ga_1)\right|+
\left|\sin({\alpha_{n_k}}\ga_2)\right|\)\[1\]\rightarrow0, \hbox{ as
} k\rightarrow\infty.\end{equation*}
In this case, we prove for
sufficiently large $k$, that \be
\frac{1}{\alpha_{n_k}}\(\frac{\partial
G}{\partial\xi}(\alpha_{n_k})\)^{-1}\leq C,~~C>0. \label{cva2g}\ee \\
Obviously, under the above assumption, we have
$|\sqrt{\mu_{n_k-2}}-\sqrt{\mu_{n_k-1}}|>\tilde{\tau}_{0}$,
$\tilde{\tau}_{0}>0$. As above, we may assume that $\mu_{n_k-1}$ and
$\mu_{n_k}$ are eigenvalues of Problems \eqref{g1} and \eqref{g2},
respectively. Hence, analogously to the first case we have the two
limits,\be
\lim_{\xi\rightarrow\sqrt{\mu_{n_k-1}}+0}F_1(\xi^2)=+\infty,~~~~
\lim_{\xi\rightarrow\sqrt{\mu_{n_k}}-0}F_2(\xi^2)=+\infty,\label{g18+*}\ee
where $F_1$ and $F_2$ are defined by \eqref{g1*}. Therefore, by
\eqref{g18+*}, for sufficiently large $k$ and each
$\alpha_{n_k}\in\(\sqrt{\mu_{n_k-1}},\sqrt{\mu_{n_k}}\),$  $$
F_1(\alpha_{n_k}^2)F_2(\alpha_{n_k}^2)>0.$$ Thus by \eqref{gap1+},
one has \begin{equation}\label{amm}\(\frac{\partial
G}{\partial\xi}(\alpha_{n_k})\)^{-1}\leq
\frac{\(F_1(\alpha_{n_k}^2)\)^2+\(F_2(\alpha_{n_k}^2)\)^2}
{2\alpha_{n_k}\(\frac{\partial
F_2}{\partial\xi}(\alpha_{n_k}^2)-\frac{\partial
F_1}{\partial\xi}(\alpha_{n_k}^2)\)}.\end{equation} In a similar way
as the first one, by use of \eqref{cvs} and \eqref{amm}, the desired
estimate
\eqref{cva2g} follows.\\
Therefore, combining \eqref{lbd}, \eqref{g20}, \eqref{g24} with
\eqref{cva2g}, we get for sufficiently large $k$,
\be\(\sqrt{{{\la}_{n_k}}}-\sqrt{{{\mu}_{n_k-1}}}\)\leq
\dfrac{C}{{n_k}},~~C>0.\label{g200}\ee
%We put\be
%\varepsilon_{n_k}^1=\sqrt{{{\la}_{n_k}}}-\sqrt{{{\mu}_{n_k-1}}}
%\hbox{ and }
%\varepsilon_{n_k}^2=\sqrt{{{\mu}_{n_k}}}-\sqrt{{{\la}_{n_k}}},\label{fg}\ee
As mentioned above, $\mu_{n_k-1}$ and $\mu_{n_k}$ are eigenvalues of
Problems \eqref{g1} and \eqref{g2}, respectively. Then by \eqref{la}
and \eqref{g21},
$$\sqrt{\la_{n_k}}\left|\tilde{u}(0,\la_{n_k})\right|\sim a_1
\left(\rho_{1}(0)\sigma_{1}(0)\right)^{-\frac{1}{4}}\gamma_1\delta_{n_k}^-,
\hbox{ as } k\rightarrow\infty,$$ where
$\delta_{n_{k}}^{-}=\sqrt{\la_{n_k}}-\sqrt{\mu_{n_{k}-1}}$. Hence,
if $\delta_{n_k}^{-}<\frac{C}{{n_k}}$, then by Lemma \ref{Lem1},\\
$\left|\tilde{u}(0,\la_{n_k})\right|\sim\left|\tilde{
v}(0,\la_{n_k})\right|$, i.e.,
$$\delta_{n_k}^{-}\sim\delta_{n_k}^{+},\hbox{ as }
k\rightarrow\infty,$$ where
$\delta_{n_{k}}^{+}=\sqrt{\mu_{n_k}}-\sqrt{\la_{n_{k}}}$. This
implies that $\left|\sqrt{\mu_{n_k}}-
\sqrt{\mu_{n_k-1}}\right|<\frac{C}{{n_k}}$, $C>0$, and this is in
contradiction with \eqref{am24}. Thus by \eqref{g200},
$$\sqrt{{{\la}_{n_k}}}-\sqrt{{{\mu}_{n_k-1}}}=\mathcal{O}(\frac{1}{n_k}),$$
which proves \eqref{g19}.\\
%\end{Proof}\\
%We are now ready to prove Theorem \ref{Tg1}.\\
%\begin{Proof}
Now we prove the estimate \eqref{g26}. If
$\frac{\gamma_{1}}{\gamma_{2}}\notin \mathbb{Q}$, then for every
$\ep>0$ and $A>0$, there exist integers $j>A$, $k>A$ such that
$$\left|\frac{j\gamma_{1}}{\gamma_{2}}-k\right|\leq \ep.$$ If
$\frac{\gamma_{1}}{\gamma_{2}}\in \mathbb{Q}$, then
$\frac{\gamma_{1}}{\gamma_{2}}=\frac{j}{k}$, for some integers
$j,~k\geq 1$. Hence, there exist two subsequences of eigenvalues
$(\sqrt{\mu^{-}_{\varphi(j)}})$ and $(\sqrt{\mu^{+}_{\varphi(k)}})$
(where $\sqrt{\mu^{-}_{n}}$ and $\sqrt{\mu^{+}_{n}}$ satisfy
\eqref{la}) such that
$$\left|\sqrt{\mu^{+}_{\varphi(k)}}-\sqrt{\mu^{-}_{\varphi(j)}}\right|
\rightarrow 0,~\hbox{as}~j,~k\rightarrow\infty,$$ or equivalently,
\begin{equation}\label{mu}
\sqrt{\mu_{\varphi(n)}}-\sqrt{\mu_{\varphi(n)-1}}\rightarrow
0,~\hbox{as}~n\rightarrow\infty.
\end{equation}
According to Proposition \ref{cg1}, there exists a subsequence
$(\la_{\varphi(n)})$ such that\\
$\la_{\varphi(n)} \in [\mu_{\varphi(n)-1}, \mu_{\varphi(n)}]$. Then
\begin{equation}\label{fff}
\sqrt{\mu_{\varphi(n)}}-\sqrt{\lambda_{\varphi(n)}}\rightarrow 0,
~\hbox{as}~n \rightarrow\infty.
\end{equation}
It is clear from \eqref{mu} and Lemma \ref{Lem1}, that
$(\varphi(n)+1)\in \Omega$, where $\Omega$ is defined by
\eqref{am4}. Hence, by \eqref{g19},
\begin{equation}\label{jh}
\sqrt{\la_{\varphi(n)+1}}-\sqrt{\mu_{\varphi(n)}}=\mathcal{O}(\frac{1}{n}).
\end{equation}
We shall prove for some $C>0$, that
\begin{equation}\label{x1}
\sqrt{\mu_{\varphi(n)}}-\sqrt{\lambda_{\varphi(n)}}\leq
\frac{C}{n},~\hbox{as}~n\rightarrow \infty.
\end{equation}
First, we suppose that $\varphi(n)\in \Omega$, then
\begin{equation}\label{ab1}
\sqrt{\la_{\varphi(n)}}-\sqrt{\mu_{\varphi(n)-1}}=
\mathcal{O}(\frac{1}{n}).
\end{equation}
As before, without loss of generality, we may assume that
$\mu_{\varphi(n)-1}$ and $\mu_{\varphi(n)}$ are eigenvalues of
Problems \eqref{g1} and \eqref{g2}, respectively. This, \eqref{la}
and \eqref{ab1}, yields
\begin{equation}\label{af1}
\left|\sin\left(\sqrt{\la_{\varphi(n)}}\gamma_{1}\right)\right|\sim\frac{C}{n},
~\hbox{for some}~C>0~\hbox{and large}~n\in \mathbb{N}^{*}.
\end{equation}
According to assertion $(ii)$ of Lemma \ref{Lem1} and \eqref{af1},
we have
\begin{equation}\label{af11}
\left|\sin\left(\sqrt{\la_{\varphi(n)}}\gamma_{2}\right)\right|\sim
\frac{C'}{n},~\hbox{for some}~C'>0.
\end{equation}
%Therefore, in view of \eqref{la} and \eqref{delta} together with
%\eqref{af1} and \eqref{af2}, one has
Since $\mu_{\varphi(n)}$ is an eigenvalue of Problem \eqref{g2},
then by \eqref{la}, \eqref{fff} and \eqref{af11},
\begin{equation}\label{af3}
\sqrt{\mu_{\varphi(n)}}-\sqrt{\lambda_{\varphi(n)}}
=\mathcal{O}(\frac{1}{n}).
\end{equation}
%Using \eqref{jh} and \eqref{af3}, the desired estimate \eqref{g26}
%holds.\\
Now, if $\varphi(n)\notin \Omega$, clearly from \eqref{la} and
\eqref{am4}, we have
$\sqrt{\mu_{\varphi(n)}}-\sqrt{\mu_{\varphi(n)-1}}< \frac{C}{n}$ for
some $C>0$ and large $n\in \mathbb{N}^{*}$. Using this and
\eqref{lam}, one gets
\begin{equation*}
\sqrt{\mu_{\varphi(n)}}-\sqrt{\la_{\varphi(n)}}<
\frac{C}{n}~\hbox{and}~\sqrt{\la_{\varphi(n)}}-\sqrt{\mu_{\varphi(n)-1}}<
\frac{C}{n}.
\end{equation*}
%Under the above assumption, since $\mu_{\varphi(n)-1}$ is an
%eigenvalue of Problem \eqref{g1}, then by \eqref{la} and
%\eqref{jhh},
%\begin{equation}\label{af11}
%\left|\sin\left(\sqrt{\la_{\varphi(n)}}\gamma_{1}\right)\right|<
%\frac{C}{n},~\hbox{for some}~C>0 .
%\end{equation}
%Hence, by Lemma \ref{Lem1},
%\begin{equation}\label{af22}
%\left|\sin\left(\sqrt{\la_{\varphi(n)}}\gamma_{2}\right)\right|<
%\frac{C}{n},~\hbox{for some}~C>0.
%\end{equation} Similarly, using
%\eqref{delta}, \eqref{af11} and \eqref{af22}, one gets
%\begin{equation}\label{af33}
%\delta_{\varphi(n)}^{-}<\frac{C}{n}~\hbox{and}~\delta_{\varphi(n)}^{+}
%<\frac{C}{n}.
%\end{equation}
Hence, from this and \eqref{af3}, the estimate \eqref{x1} holds.
Combining \eqref{jh} and \eqref{x1}, the desired
estimate \eqref{g26} follows.\\
It is easily seen from \eqref{lam} and assertion $(i)$ of Lemma
\ref{Lem1}, that the estimate \eqref{g27} holds. This ends up the
proof of the theorem.
\end{Proof}\\
Proceeding as above, we have the following result:
\begin{coro}\label{cor3}
Let $(\mu_{n_k})$ be a subsequence of $(\mu_{n})$.\\
If $\sqrt{\mu_{n_{k}}}-\sqrt{\mu_{n_{k}-1}}\rightarrow 0$ as $k
\rightarrow \infty$, then
$\sqrt{\la_{n_{k}+1}}-\sqrt{\la_{n_{k}}}=\mathcal{O}(\frac{1}{k})$,
$\sqrt{\mu_{n_k}}-\sqrt{\la_{n_{k}}}\leq \frac{C}{k}$ and
$\sqrt{\la_{n_{k}}}-\sqrt{\mu_{n_{k}-1}}\leq \frac{C}{k}$ for some
$C>0$.\\
Moreover, if $|\sqrt{\mu_{n_{k}}}-\sqrt{\mu_{n_{k}-1}}|> C$ for some
$C>0$ and large $k$,
then\\$|\sqrt{\la_{n_{k}+1}}-\sqrt{\la_{n_{k}}}|
> C'$ for all $k$.
\end{coro}
We establish now the asymptotic behavior of the eigenfunctions
$(\phi_{n}(x))_{n\geq1}$ of the eigenvalue problem \eqref{an}.
\begin{prop}\label{ff}
Define the set
\begin{equation}\label{k}
\Lambda=\left\{n\in\mathbb{N^{*}}~\hbox{such
that}~\mu_{n}\in\Gamma^{*} \right\},\end{equation} where the set
$\Gamma^{*}$ is defined by \eqref{g34}. Then the associated
eigenfunctions $(\phi_{n}(x))_{n\geq1}$ of the eigenvalue problem
\eqref{an} satisfy the following
asymptotic estimates:\\
\begin{description}
  \item[i$)$] For $n\in \Lambda$,
\bean\label{pro1}\phi_{n}(x)= \left\{
\begin{array}{ll}
-a_{1}a_{2}\left(\dfrac{\rho_{1}(x)\sigma_{1}(x)}
{\rho_{2}(0)\sigma_{2}(0)}\right)^{-\frac{1}{4}}\dfrac{\cos(\sqrt{
\lambda_{n}}\gamma_{2})\sin\left(\sqrt{\lambda_{n}}\int_{-1}^{x}
\sqrt{\frac{\rho_{1}(t)}{\sigma_{1}(t)}}dt\right)}
{\sqrt{\lambda_{n}}}\[1\],~~x\in [-1,0],\\
a_{1}a_{2}\left(\dfrac{\rho_{2}(x)\sigma_{2}(x)}{\rho_{1}(0)\sigma_{1}(0)}\right)^{-\frac{1}{4}}\dfrac{\cos(\sqrt{\lambda_{n}}\gamma_{1})\sin\left(\sqrt{\lambda_{n}}\int_{x}^{1}
\sqrt{\frac{\rho_{2}(t)}{\sigma_{2}(t)}}dt\right)}{\sqrt{\lambda_{n}}
}\[1\],~~x\in [-1,0],
\end{array}
 \right.
\eean where $\[1\]=1+\mathcal{O}(\frac{1}{n})$,
$\gamma_{1}=\int_{-1}^{0}\sqrt{\frac{\rho_{1}(x)}{\sigma_{1}(x)}}dx,$
 $\gamma_{2}=\int_{0}^{1}\sqrt{\frac{\rho_{2}(x)}{\sigma_{2}(x)}}dx$
 and the constants $a_{i}~~(i=1,2)$
are defined by \eqref{am27}.
\item[ii$)$] For $n \in \mathbb{N}^{*}\setminus\Lambda$,
   \bean\label{pro2} \phi_{n}(x)=\left\{
\begin{array}{ll}
a_{1}a_{2}\dfrac{\left(\rho_{1}(x)\sigma_{1}(x)\right)^{-\frac{1}{4}}}{\left(\rho_{2}(0)\sigma_{2}(0)\right)^{\frac{1}{4}}}
\dfrac{\sin(\sqrt{\lambda_{n}}\gamma_{2})\sin\left(\sqrt{\lambda_{n}}\int_{-1}^{x}
\sqrt{\frac{\rho_{1}(t)}{\sigma_{1}(t)}}dt\right)}{\sqrt{\lambda_{n}}}
\[1\],~~x\in [0,1],\\
a_{1}a_{2}\dfrac{\left(\rho_{2}(x)\sigma_{2}(x)\right)^{-\frac{1}{4}}}{\left(\rho_{1}(0)\sigma_{1}(0)\right)^{\frac{1}{4}}}\dfrac{\sin(\sqrt{\lambda_{n}}\gamma_{1})
\sin\left(\sqrt{\lambda_{n}} \int_{x}^{1}
\sqrt{\frac{\rho_{2}(t)}{\sigma_{2}(t)}}dt\right)}
{\sqrt{\lambda_{n}}}\[1\],~~x\in [0,1].
\end{array}
 \right.
\eean
\end{description}
\end{prop}
\begin{Proof}
It is clear that for all $n\in\Lambda$,
$\tilde{u}(0,\lambda_{n})=\tilde{v}(0,\lambda_{n})=0$. Then, from
\eqref{ll} together with the last condition in \eqref{an}, the
corresponding eigenfunctions $(\phi_{n}(x))_{n\in\Lambda}$ can be
written in the form
\begin{equation}\label{ss}
\phi_{n}(x)= \left\{
\begin{array}{lll}
\sigma_{2}(0)\tilde{v}_{x}(0,\lambda_{n})\tilde{u}(x,\lambda_{n}),&-1\leq x\leq 0,\\
\sigma_{1}(0)\tilde{u}_{x}(0,\lambda_{n})\tilde{v}(x,\lambda_{n}),&0\leq x\leq 1.\\
\end{array}
\right.
\end{equation}
Therefore, the asymptotes \eqref{pro1} are simple deductions from
the asymptotes \eqref{g21}, \eqref{g22} and \eqref{ss}. Now, if
$n\in\mathbb{N^{*}}\backslash\Lambda$, we have
$\tilde{u}(0,\lambda_{n})\neq 0$ and $\tilde{v}(0,\lambda_{n})\neq
0$, since otherwise $\la_n$ is not an eigenvalue of Problem
\eqref{an} (see the second case in the proof of Lemma \ref{lemm2}).
We set
\begin{equation}\label{gg}
\phi_{n}(x)=\sqrt{\lambda_{n}}\widetilde{U}(x,\la_{n}),
\end{equation}
where $\widetilde{U}(x,\la_{n})$ is defined by \eqref{g4}. From the
asymptotes \eqref{g21}, \eqref{g22} and \eqref{gg}, a
straightforward computation gives the asymptotes \eqref{pro2}.
\end{Proof}
\section{Riesz Basis And Asymmetric Space} In this section we give
results concerning the Riesz basis and the asymmetric space that we
will need later to prove the observability inequality. Let us begin
with the following result due to C. Baiocchi et al. \cite{Bi1}:
\begin{theo}\label{cor1}
Let $(\omega_{n})$ be a strictly increasing sequence satisfying for
some $\delta >0$ the condition:
$$\omega_{n+2}-\omega_{n} > 2\delta,~\hbox{for all}~n.$$
Fix a number $0<\delta^{'}\leq \delta$, and set
$$A=\left\{n\in \mathbb{Z},~~\omega_{n}-\omega_{n-1}\geq \delta^{'}~
\hbox{and}~ \omega_{n+1}-\omega_{n}<\delta^{'}\right\},$$
$$B=\left\{n\in \mathbb{Z},~~ \omega_{n}-\omega_{n-1}\geq \delta^{'}~
\hbox{and}~ \omega_{n+1}-\omega_{n}\geq \delta^{'}\right\}.$$ Then
the following estimates hold for every bounded interval I of length
$|I|> 2\pi D^{+}$: there exist two constants $C_{1}>0$ and $C_{2}>0$
such that for any $f=\sum\limits_{n \in \mathbb{Z}} a_{n} e^{i
\omega_{n}t}$,
\begin{align*}
&C_{1} \left(\sum_{n \in A}
\left[(\omega_{n+1}-\omega_{n})^{2}(|a_{n}|^{2}+|a_{n+1}|^{2})
+|a_{n+1}+a_{n}|^{2}\right]
+\sum_{n \in B}|a_{n}|^{2}\right)\leq \int_{I}|f(t)|^{2}dt\\
&\leq C_{2} \left( \sum_{n \in A}
\left[(\omega_{n+1}-\omega_{n})^{2}(|a_{n}|^{2}+|a_{n+1}|^{2})+|a_{n+1}+a_{n}|^{2}\right]
+\sum_{n \in B}|a_{n}|^{2} \right),
\end{align*}
where $D^{+}:=\lim\limits_{r\rightarrow +\infty} \frac{n^{+}(r)}{r}$
and $n^{+}(r)$ denotes the largest number of terms of the sequence
$(\omega_{n})$ contained in an interval of length $r$.
\end{theo}
%From Theorem \ref{Tg1}, there exists a constant $\delta >0$ such
%that $\sqrt{\lambda_{n+2}}-\sqrt{\lambda_{n}} \geq 2\delta$, for all
%$n\in \mathbb{N}^{*}$. Moreover, there exists a subsequence of
%eigenvalues $(\lambda_{\varphi(n)})$ such that
%$\sqrt{\lambda_{\varphi(n)+1}}-\sqrt{\lambda_{\varphi(n)}}=\mathcal{O}(\frac{1}{n})$.
Set
\begin{equation}\label{zed}
\sqrt{\lambda_{-n}}=-\sqrt{\lambda_{n}},~\sqrt{\mu_{-n}}=-\sqrt{\mu_{n}}
~\hbox{and}~\Phi_{-n}=\Phi_{n},~ \hbox{for all}~n\in \mathbb{N}^{*}.
\end{equation}
In light of Theorem \ref{Tg1} and Corollary \ref{cor3}, let
$\sigma^{*}$ be the set of elements of all subsequences
$(\la_{n_k})$ which satisfy
$\sqrt{\la_{n_{k}+1}}-\sqrt{\la_{n_k}}=\mathcal{O}(\frac{1}{k})$,
i.e.,
\begin{equation*}
\sigma^{*}=\{\la_{n_k}~\hbox{such
that}~\sqrt{\la_{n_{k}+1}}-\sqrt{\la_{n_k}}=\mathcal{O}(\frac{1}{k})\}.
\end{equation*}
Denote by
\begin{equation}\label{aa1}
A=\left\{n\in \mathbb{Z}^{*}~\hbox{such that}~\la_{n}\in
\sigma^{*}\right\}
\end{equation}and
\begin{equation}\label{aa2}
B=\left\{n\in \mathbb{Z}^{*},~(n-1)\notin A~\hbox{and}~n\notin
A\right\}.
\end{equation}
Note that if $n\in A$, then $(n+1)\not\in A$ and $(n+1)\not\in B$.
Thus
\begin{equation}\label{hn}
\mathbb{Z}^{*}=A\cup B \cup\{n+1,~~n\in A\}.
\end{equation}
Moreover, it is easily seen from Corollary \ref{cor3}, that for
large $n\in A$, $|\sqrt{\mu_{n}}-\sqrt{\mu_{n-1}}|$ tends to zero.
%\begin{rem}\label{rem1}
%If there exists an other subsequence $(\sqrt{\mu_{\psi(n)}})$ such
%that\\$\sqrt{\mu_{\psi(n)}}-\sqrt{\mu_{\psi(n)-1}}\rightarrow 0$ as
%$n \rightarrow \infty$, then proceeding as above, we prove that\\
%$\sqrt{\la_{\psi(n)+1}}-\sqrt{\la_{\psi(n)}}=\mathcal{O}(\frac{1}{n})$.
%Hence for $n\in A$, we have
%$\sqrt{\la_{n+1}}-\sqrt{\la_{n}}=\mathcal{O}(\frac{1}{n})$.
%\end{rem}
\begin{lemm}\label{cor2}
For $T> 2(\gamma_{1}+\gamma_{2})$, where
$\gamma_{1}=\int_{-1}^{0}\sqrt{\frac{\rho_{1}(x)}{\sigma_{1}(x)}}dx$
and
$\gamma_{2}=\int_{0}^{1}\sqrt{\frac{\rho_{2}(x)}{\sigma_{2}(x)}}dx$,
there exist two positive constants $C_{1}$ and $C_{2}$ such that for
$f=\sum\limits_{n\in \mathbb{Z}^{*}}a_{n}e^{i\sqrt{\lambda_{n}}t}$,
we have
\begin{align}
&C_{1} \left(\sum_{n \in A}
\left[(|a_{n}|^{2}+|a_{n+1}|^{2})\delta_{n}^{2}+|a_{n+1}+a_{n}|^{2}\right]
+\sum_{n \in B}|a_{n}|^{2}\right)\leq \int_{0}^{T}\left|\sum_{n\in
\mathbb{Z}^{*}}a_{n}e^{i\sqrt{\lambda_{n}}t}\right|^{2}dt\nonumber\\
&\leq C_{2} \left(\sum_{n \in A}
\left[(|a_{n}|^{2}+|a_{n+1}|^{2})\delta_{n}^{2}+|a_{n+1}+a_{n}|^{2}\right]
+\sum_{n \in B}|a_{n}|^{2} \right).\label{hhh}
\end{align}
\end{lemm}
\begin{Proof}  By \eqref{g27}, we have for $n\in A$,
$\sqrt{\la_{n}}-\sqrt{\la_{n-1}}\geq\delta'$ for some $\delta'<
\delta$. Hence the second condition in the set $A$ given in Theorem
\ref{cor1} holds for $\omega_{n}=\sqrt{\la_{n}}$.\\
In view of Proposition \ref{cg1},
$\frac{n^{+}(\sqrt{\lambda^{'}_{n-1}})}{\sqrt{\lambda^{'}_{n}}}\leq
\frac{n^{+}(\sqrt{\lambda_{n}})}{\sqrt{\lambda_{n}}} \leq
\frac{n^{+}(\sqrt{\lambda^{'}_{n}})}{\sqrt{\lambda^{'}_{n-1}}},$
where $\lambda'_{n}$ $(n\geq1)$ denote the eigenvalues of the
regular problem \eqref{an} for $M=0$, and $n^{+}$ is defined in
Theorem \ref{cor1}. Using the asymptote \eqref{g23}, we find $
D^{+}= \frac{\gamma_{1}+\gamma_{2}}{\pi}$. By setting
$\omega_{n}=\sqrt{\lambda_{n}}$,\\
$\delta_{n}=\sqrt{\lambda_{n+1}}-\sqrt{\lambda_{n}}$, $n\in
\mathbb{Z}^{*}$, in view of \eqref{g26} and \eqref{g27}, we are in
the conditions of Theorem \ref{cor1}, and hence, Inequality
\eqref{hhh} follows.
\end{Proof}\\
Define the scale of Hilbert spaces $(X_{\alpha}, \|.\|_{\alpha})$,
$\alpha \in \mathbb{R}$:
\begin{equation}\label{xa}
X_{\alpha}=\left\{Y: Y=\sum_{n\in
\mathbb{N}^{*}}a_{n}\Phi_{n}~\hbox{with}~\|Y\|_{\alpha}^{2}=
\sum_{n\in\mathbb{N}^{*}}|a_{n}|^{2}\lambda^{2\alpha}_{n}\langle
\Phi_{n},\Phi_{n}\rangle_{\mathcal{W}} < \infty\right\}
\end{equation}endowed with the norm $\|.\|_{\alpha}$, where
$\mathcal{W}$ is defined by \eqref{v} and
$\Phi_{n}$ are the eigenfunctions of the operator $\mathcal{A}$.\\
\begin{prop}\label{prop11} We have
\begin{equation}\label{k1}
\langle \Phi_{n},\Phi_{n}\rangle_{\mathcal{W}}= \mathcal{O}(1),~n\in
\Lambda\cup (B\cap \mathbb{N}^{*}),
\end{equation}
\begin{equation}\label{k2}
\langle \Phi_{n},\Phi_{n}\rangle_{\mathcal{W}}\leq
\frac{C}{n^{2}},~\hbox{for large}~n\in (A\cap
\mathbb{N}^{*})\setminus \Lambda~\hbox{and}~C>0,
\end{equation}
and
\begin{equation}\label{k3}
\langle \Phi_{n+1},\Phi_{n+1}\rangle_{\mathcal{W}}
=\mathcal{O}(\frac{1}{n^{2}}),~\hbox{for large}~n\in A\cap
\mathbb{N}^{*}.
\end{equation}
where the sets $\Lambda$, $A$ and $B$ are defined by \eqref{k},
\eqref{aa1} and \eqref{aa2}, respectively.
\end{prop}
\begin{Proof}
In view of \eqref{j11}, \eqref{g21} and \eqref{g22}, we have for
$n\in \mathbb{N}^{*}$
\begin{equation}\label{k5}
\|\tilde{u}(x,\la_{n})\|^{2}_{\mathcal{V}_{1}}=\mathcal{O}(1)~\hbox{and}~
\|\tilde{v}(x,\la_{n})\|^{2}_{\mathcal{V}_{2}}=\mathcal{O}(1).
\end{equation}
It is clear from \eqref{j33} and \eqref{ss}, that for $n\in \Lambda$
\begin{equation}\label{k4}
\langle \Phi_{n},\Phi_{n}\rangle_{\mathcal{W}}=
\left(\sigma_{2}(0)\tilde{v}_{x}(0,\lambda_{n})\right)^{2}
\|\tilde{u}(x,\la_{n})\|^{2}_{\mathcal{V}_{1}}+
\left(\sigma_{1}(0)\tilde{u}_{x}(0,\lambda_{n})\right)^{2}
\|\tilde{v}(x,\la_{n})\|^{2}_{\mathcal{V}_{2}}.
\end{equation}
Note that for $n\in \Lambda$, $\la_{n}=\mu_{n}$, i.e.,
$\tilde{u}(0,\lambda_{n})=\tilde{v}(0,\lambda_{n})=0$. In this case
$\la_n$ is an eigenvalue of both Problems \eqref{g1} and \eqref{g2}.
This means that $\tilde{u}_{x}(0,\lambda_{n})\neq0$ and
$\tilde{v}_{x}(0,\lambda_{n})\neq0$ for all $n\in \Lambda$, and
hence, by \eqref{la}, \eqref{g21} and \eqref{g22}, we get for large
$n\in \Lambda$
\begin{equation}\label{sig1}
\left|\sigma_{1}(0)\tilde{ u}_{x}(0,\la_{n})\right|\sim
\frac{\left(\rho_{1}(0)\sigma_{1}(0)\right)^{\frac{1}{4}}}
{\left(\rho_{1}(-1)\right)^{\frac{1}{4}}
\left(\sigma_{1}(-1)\right)^{-\frac{3}{4}}}~\hbox{and}~
\left|\sigma_{2}(0)\tilde{v}_{x}(0,\la_{n})\right|\sim
\frac{\left(\rho_{2}(0)\sigma_{2}(0)\right)^{\frac{1}{4}}}
{\left(\rho_{2}(1)\right)^{\frac{1}{4}}
\left(\sigma_{2}(1)\right)^{-\frac{3}{4}}}.\end{equation} From this,
\eqref{k5} and
\eqref{k4}, the estimate \eqref{k1} holds for $n\in \Lambda$.\\
Now for $n\in \mathbb{N}^{*} \setminus \Lambda$, using \eqref{j33}
and \eqref{gg}, we obtain
\begin{equation}\label{k44}
\langle \Phi_{n},\Phi_{n}\rangle_{\mathcal{W}}=
\left(\sqrt{\la_{n}}\tilde{v}(0,\lambda_{n})\right)^{2}
\|\tilde{u}(x,\la_{n})\|^{2}_{\mathcal{V}_{1}}+
\left(\sqrt{\la_{n}}\tilde{u}(0,\lambda_{n})\right)^{2}
\|\tilde{v}(x,\la_{n})\|^{2}_{\mathcal{V}_{2}}.
\end{equation}
By Corollary \ref{cor3}, we have for large $n\in (A\cap
\mathbb{N}^{*})\setminus \Lambda$,
\begin{equation}\label{b2}
|\sqrt{\la_{n}}-\sqrt{\mu_{n}}|\leq \frac{C}{n}~\hbox{and}~
|\sqrt{\la_{n}}-\sqrt{\mu_{n-1}}|\leq \frac{C}{n},~\hbox{for
some}~C>0.
\end{equation}
Then by use of \eqref{la}, \eqref{g21} and \eqref{g22} together with
\eqref{b2}, we obtain
\begin{equation}\label{v0}
\left|\sqrt{\lambda_{n}}\tilde{u}(0,\lambda_{n})\right|\leq
\frac{C'}{n}~\hbox{and}~
\left|\sqrt{\lambda_{n}}\tilde{v}(0,\lambda_{n})\right|\leq
\frac{C''}{n},~\hbox{for some}~C',~C''>0.
\end{equation}
Since for large $n\in (A \cap \mathbb{N}^{*})$,
$|\sqrt{\mu_{n}}-\sqrt{\mu_{n-1}}|$ tends to zero, then in view of
Corollary \ref{cor3} together with \eqref{g19},
\begin{equation*}
\sqrt{\la_{n+1}}-\sqrt{\mu_{n}}=\mathcal{O}(\frac{1}{n})~\hbox{and}~
\sqrt{\la_{n+1}}-\sqrt{\mu_{n-1}}=\mathcal{O}(\frac{1}{n}),~\hbox{for
large}~ n\in (A\cap \mathbb{N}^{*}).
\end{equation*}
Substituting these asymptotes into \eqref{g21} and \eqref{g22}, one
gets
\begin{equation}
\sqrt{\lambda_{n+1}}\tilde{u}(0,\lambda_{n+1})=
\mathcal{O}(\frac{1}{n})~\hbox{and}~\sqrt{\lambda_{n+1}}\tilde{v}(0,\lambda_{n+1})=
\mathcal{O}(\frac{1}{n}).\label{vv0}
\end{equation}
It is easily seen that $|\sqrt{\mu_{n}}-\sqrt{\mu_{n-1}}|>\delta'$
for all $n\in B$. Then again by \eqref{la} and \eqref{g21} together
with \eqref{g22} and \eqref{g19},
\begin{equation*}
\left|\sqrt{\lambda_{n}}\tilde{u}(0,\lambda_{n})\right|+
\left|\sqrt{\lambda_{n}}\tilde{v}(0,\lambda_{n})\right|>
C,~\hbox{for some}~C>0~\hbox{and large}~n\in B.
\end{equation*}
Therefore, from the above and \eqref{k44}, the desired estimates
\eqref{k1}, \eqref{k2} and \eqref{k3} hold. The proposition is
proved.
\end{Proof}\\
As a consequence of Proposition \ref{prop11}, we give other
descriptions
of some of the spaces $X_{\alpha}$.%According to Proposition \ref{33},
%$X_{0}$ coincides with the subspace $\mathcal{V}$, where
%$\mathcal{V}$ is defined by (\ref{w3}). On the other hand,
%$X_{-1/2}=\mathcal{H}$ and $X_{-1}$ coincides with the dual space of
%$D(\mathcal{A}^{\frac{1}{2}})$. Thus
%$$X_{-1}=H^{-1}(-1,0)\times H^{-1}(0,1)\times \mathbb{R}.$$
%In view of Proposition \ref{33}, $\|.\|_{\mathcal{W}}$ is equivalent
%to $\|.\|_{\mathcal{V}}$ which is defined by (\ref{n1}).
\begin{prop}\label{pr}
We have the following characterizations of the spaces $X_{\alpha}$:
\begin{enumerate}
\item $X_{0}=D(\mathcal{A}^{\frac{1}{2}})$ coincides topologically with the
subspace $\mathcal{W}$, where $\mathcal{W}$ is defined by
\eqref{v}.\\
\item $X_{-\frac{1}{2}}=D(\mathcal{A}^{0})$ coincides algebraically and topologically with the
space $\mathcal{H}_{0}$, where $\mathcal{H}_{0}$ is defined by \eqref{h}.\\
\item $X_{-1}=D(\mathcal{A}^{-\frac{1}{2}})$
coincides with the dual space of $X_{0}$, it is the subspace
$\mathcal{H}_{-1}$, where $\mathcal{H}_{-1}$ is defined by
\eqref{h1}.
\end{enumerate}
\end{prop}
\begin{Proof}
The proof is analogous to that of \cite{Beam2}.
\end{Proof}\\
%Throughout the rest of the paper we denote by
%$\overline{f}(x,\lambda)=(f, i \sqrt{\lambda}f)$.
Let us now recall briefly how solutions of
\eqref{a}-\eqref{22}-\eqref{j} can be developed in Fourier series.
Given initial data $(u^{0},v^{0},z^{0})\in X_{-\frac{1}{2}}$,
$(u^{1},v^{1},z^{1})\in X_{-1}$, we compute its Fourier coefficients
\begin{equation}\label{xx1}
\tilde{e}_{n}=\frac{\langle
(u^{0},v^{0},z^{0}),\Phi_{n}\rangle_{X_{-\frac{1}{2}}}}{\langle
\Phi_n, \Phi_n \rangle_{X_{-\frac{1}{2}}}},~~
\tilde{f}_{n}=\frac{\langle
(u^{1},v^{1},z^{1}),\Phi_{n}\rangle_{X_{-1}}}{\langle \Phi_n, \Phi_n
\rangle_{X_{-\frac{1}{2}}}},
\end{equation}
where $\Phi_{n}=(\phi_{n},\phi_{n}(0))$ and $\langle \Phi_{n},
\Phi_n\rangle_{X_{-\frac{1}{2}}}=\frac{\langle \Phi_{n},
\Phi_n\rangle_{
\mathcal{W}}}{\la_n}$. %Thus
%\begin{equation*}
%u^{0}(x)=\sum_{n\in \mathbb{N}^{*}}\tilde{e}_n\phi_{n}^{u},
%~v^{0}(x)=\sum_{n\in
%\mathbb{N}^{*}}\tilde{e}_n\phi_{n}^{v},~z^{0}=\sum_{n\in
%\mathbb{N}^{*}}\tilde{e}_n\phi_{n}(0),
%\end{equation*}
%\begin{equation*}
%u^{1}(x)=\sum_{n\in
%\mathbb{N}^{*}}\tilde{f}_n\sqrt{\la_n}\phi_{n}^{u},~
%v^{1}(x)=\sum_{n\in
%\mathbb{N}^{*}}\tilde{f}_n\sqrt{\la_n}\phi_{n}^{v},~z^{1}=\sum_{n\in
%\mathbb{N}^{*}}\tilde{f}_n\sqrt{\la_n}\phi_{n}(0),
%\end{equation*}
%where $\phi_{n}^{u}$ and $\phi_{n}^{v}$ denote the restrictions of
%the eigenfunctions $\phi_{n}$ to $(-1,0)$ and $(0,1)$, respectively.
%Then the unique solution $(u,v,z)$ of \eqref{a}-\eqref{22}-\eqref{j}
%can be written as
%\begin{equation*}
%u(x,t)=\sum_{n\in
%\mathbb{N}^{*}}\left(\tilde{e}_n\cos(\sqrt{\la_n}t)\phi_{n}^{u}(x)+
%\tilde{f}_n\sin(\sqrt{\la_n}t)\phi_{n}^{u}(x)\right),
%\end{equation*}
%\begin{equation*}
%v(x,t)=\sum_{n\in
%\mathbb{N}^{*}}\left(\tilde{e}_n\cos(\sqrt{\la_n}t)\phi_{n}^v(x)+
%\tilde{f}_n\sin(\sqrt{\la_n}t)\phi_{n}^{v}(x)\right),
%\end{equation*}
%\begin{equation*}
%z(t)=\sum_{n\in
%\mathbb{N}^{*}}\left(\tilde{e}_n\cos(\sqrt{\la_n}t)\phi_{n}(0)+
%\tilde{f}_n\sin(\sqrt{\la_n}t)\phi_{n}(0)\right).
%\end{equation*}
%To simplify these expressions,
Now, we introduce the complex Fourier coefficients
\begin{equation}\label{xx2}
a_{n}=\frac{\tilde{e}_n-i\tilde{f}_n}{2},~~
a_{-n}=\frac{\tilde{e}_n+i\tilde{f}_n}{2},~n\in
\mathbb{N}^{*}.\end{equation} Whence, from this and the notations of
\eqref{zed}, the solution $(u,v,z)$ can be written as follows
$$(u,v,z)=\sum_{n\in
\mathbb{Z}^{*}}a_{n}e^{i\sqrt{\la_{n}}t}\Phi_{n}.$$ According to
Lemma \ref{lemm2}, it is clear that the set
$(\overline{\Phi}_{n})_{n\in \mathbb{Z}^{*}}$, with
$\overline{\Phi}_n=(\Phi_n,i\sqrt{\la_n}\Phi_n)$, forms an
orthogonal basis of the energy space $X_{\alpha}\times
X_{\alpha-\frac{1}{2}}$. Then from the above, the vector valued
solution $U=((u,v,z),(u_{t},v_{t},z_{t}))$ of System
\eqref{a}-\eqref{22}-\eqref{j} is given by
$$U(t)=\sum_{n\in
\mathbb{Z}^{*}}a_{n}e^{i\sqrt{\lambda_{n}}t}\overline{\Phi}_{n}.$$
%where the complex coefficient $a_{n}$ is determined by the initial
%data \eqref{j}.\\ %In the expression above $\Phi_{n}$ is identified
%with the vector
%$(\phi_{n}^{u}(x),\phi_{n}^{v}(x),\phi_{n}^{u}(0))$.\\
Now we introduce a subspace of the energy space $X_{-\frac{1}{2}}
\times X_{-1}$ defined by:
\begin{align}
\mathcal{Y}&=\big\{U=\sum_{n\in \mathbb{Z}^{*}} a_{n}
\overline{\Phi}_{n} \in X_{-\frac{1}{2}} \times
X_{-1}:~\|U\|^{2}_{\mathcal{Y}}=\sum_{n\in
A}\delta_{n}^{2}\left(|\tilde{a}_{n}|^{2}+|\tilde{a}_{n+1}|^{2}
\right)+|\tilde{a}_{n}+\tilde{a}_{n+1}|^{2}\nonumber\\
&+\sum_{n\in B}|\tilde{b}_{n}|^{2}< \infty \big\},\label{y}
\end{align}
where $\tilde{a}_{n}=\phi_{n}'(1)a_{n}$ and
$\tilde{b}_{n}=\phi_{n}'(1)b_{n}$. Clearly by Proposition \ref{pr},
$\mathcal{Y}$ is a subspace of the energy space
$X_{-\frac{1}{2}}\times X_{-1}=\mathcal{H}_{0}\times
\mathcal{H}_{-1}$, where $\mathcal{H}_{0}$ and $\mathcal{H}_{-1}$
are defined by \eqref{h} and \eqref{h1}, respectively. Let us see
that System \eqref{a}-\eqref{22}-\eqref{j} is well-posed in
$\mathcal{Y}$.
\begin{lemm}
Let $U^{0}=((u^{0},v^{0},z^{0}),(u^{1},v^{1},z^{1}))$ be an element
of $\mathcal{Y}$. Then the solution $U$ of \eqref{a}-\eqref{22} with
initial data $U^{0}$ belongs to $\mathcal{Y}$ for every $t>0$.
Furthermore, for any $T>0$ there exists a constant $C(T)>0$ such
that
$$\|U(t)\|_{\mathcal{Y}} \leq C(T)
\|U^{0}\|_{\mathcal{Y}},~~~~~~0\leq t\leq T~and~ U^{0}\in
\mathcal{Y}.$$
\end{lemm}
\begin{Proof}
Given $U^{0}\in \mathcal{Y}$, the unique solution $U$ of
\eqref{a}-\eqref{22} with initial data $U^{0}$ can be represented in
Fourier series as follows
$$U=\sum\limits_{n\in \mathbb{Z}^{*}} a_{n}\overline{\Phi}_{n}(x)
e^{i\sqrt{\lambda_{n}}t}.$$ The Fourier coefficients $a_{n}$ are
determined by the initial data $U^{0}=\sum\limits_{n\in
\mathbb{Z}^{*}} a_{n}\overline{\Phi}_{n}(x)$. On the other hand,
\begin{align*}
&\|U(t)\|^{2}_{\mathcal{Y}}=\sum_{n\in
A}\delta_{n}^{2}\left(|\tilde{a}_{n}e^{i\sqrt{\lambda_{n}}t}|^{2}
+|\tilde{a}_{n+1}e^{i\sqrt{\lambda_{n+1}}t}|^{2}\right)
+|\tilde{a}_{n}e^{i\sqrt{\lambda_{n}}t}+\tilde{a}_{n+1}e^{i\sqrt{\lambda_{n+1}}t}|^{2}\\
&+\sum_{n\in B}|\tilde{b}_{n}e^{i\sqrt{\lambda_{n}}t}|^{2}\leq 2
\sum_{n\in A}\delta_{n}^{2}\left(|\tilde{a}_{n}|^{2}+
|\tilde{a}_{n+1}|^{2}\right)+|\tilde{a}_{n+1}|^{2}
|e^{i\sqrt{\lambda_{n+1}}t}-e^{i\sqrt{\lambda_{n}}t}|^{2}\\
&+2 \sum_{n\in A}|\tilde{a}_{n}+
\tilde{a}_{n+1}|^{2}|e^{i\sqrt{\lambda_{n}}t}|^{2} +\sum_{n\in
B}|\tilde{b}_{n}|^{2}.
\end{align*}
Since for $n\in A$,
$\big|e^{i\sqrt{\lambda_{n+1}}t}-e^{i\sqrt{\lambda_{n}}t}\big|=(\sqrt{\la_{n+1}}-\sqrt{\la_n})t
=\delta_{n}t$, then
\begin{align*}
&\|U(t)\|^{2}_{\mathcal{Y}}\leq C \left(\sum_{n\in
A}\delta_{n}^{2}\left(|\tilde{a}_{n}|^{2}+(1+t^{2})
|\tilde{a}_{n+1}|^{2}\right)+|\tilde{a}_{n}+\tilde{a}_{n+1}|^{2}\right)\\
&+C\sum_{n\in B}|\tilde{b}_{n}|^{2}\leq
C(T)\|U^{0}\|^{2}_{\mathcal{Y}}.
\end{align*}
\end{Proof}\\
In order to characterize the space $\mathcal{Y}$, we introduce for
$n\in A$
%\begin{equation}\label{phi}
%\tilde{\Phi}_{n}=\frac{\Phi_{n}} {\phi_{n}'(1)},
%\end{equation}
\begin{equation}\label{cc1}
q_{n}=\frac{1}{2}\left(
\frac{\overline{\Phi}_{n+1}}{\phi_{n+1}'(1)}+
\frac{\overline{\Phi}_{n}}{\phi_{n}'(1)}\right)
~\hbox{and}~p_{n}=\frac{1}{2\delta_{n}}\left(
\frac{\overline{\Phi}_{n+1}}{\phi_{n+1}'(1)}-
\frac{\overline{\Phi}_{n}}{\phi_{n}'(1)}\right).
\end{equation}
\begin{prop}\label{44}
The set $(p_{n})_{n\in A}\cup (q_{n})_{n\in A} \cup
\left(\frac{\overline{\Phi}_{n}}{\phi_{n}'(1)}\right)_{n\in B}$
forms a Riesz basis in the space $\mathcal{Y}$.
\end{prop}
\begin{Proof}
From \eqref{hn}, observe that any element $U\in \mathcal{Y}$ can be
written as
\begin{align}
&U=\sum_{n\in \mathbb{Z}^{*}}a_{n}\overline{\Phi}_{n}=\sum_{n\in
A}a_{n}\overline{\Phi}_{n}+a_{n+1}\overline{\Phi}_{n+1}+\sum_{n\in
B}b_{n}\overline{\Phi}_{n}\nonumber\\
&=\sum_{n\in A}(\tilde{a}_{n}+\tilde{a}_{n+1})q_{n}+\delta_{n}
(\tilde{a}_{n+1}-\tilde{a}_{n})p_{n}+\sum_{n\in
B}\tilde{b}_{n}\frac{\overline{\Phi}_{n}}{\phi_{n}'(1)},\label{ghh}
\end{align}
where $\tilde{a}_{n}=a_{n}\phi_{n}'(1)$ and
$\tilde{b}_{n}=b_{n}\phi_{n}'(1)$.
Whence, the set under
consideration is complete in $\mathcal{Y}$. Moreover, if we define
on $\mathcal{Y}$ a scalar product such that this set is orthonormal,
then the corresponding norm is such that
\begin{align*}
&\|U\|_{*}^{2}=\sum_{n\in
A}\delta_{n}^{2}|\tilde{a}_{n+1}-\tilde{a}_{n}|^{2}
+|\tilde{a}_{n+1}+\tilde{a}_{n}|^{2} +\sum_{n\in
B}|\tilde{b}_{n}|^{2}.
\end{align*}
%\begin{align*}
%&\|U\|^{2}_{\mathcal{Y}}\leq C_{1}\sum_{n\in
%\Lambda^{*}}\delta_{n}^{2}\left(|a_{n+1}\phi_{n+1}'(1)|^{2}
%+|a_{n}\phi_{n}'(1)|^{2}\right)+|a_{n+1}\phi_{n+1}'(1)+a_{n}\phi_{n}'(1)|^{2}\\
%&+C_{2}\sum_{n\in A \setminus
%\Lambda^{*}}\delta_{n}^{2}\left(|a_{n+1}\phi_{n+1}'(1)|^{2}
%+|a_{n}\phi_{n}'(1)|^{2}\right)
%+|a_{n+1}\phi_{n+1}'(1)+a_{n}\phi_{n}'(1)|^{2}\\
%&+C_{3}\sum_{n\in B^{+}}|b_{n}\phi_{n}'(1)|^{2}+C_{4}\sum_{n\in
%B^{-}}|b_{n}\phi_{n}'(1)|^{2}.
%\end{align*}
Since $\delta_{n}=\mathcal{O}(\frac{1}{n})$, it is easily seen that
\begin{align*}
&\|U\|^{2}_{\mathcal{Y}}\leq \sum_{n\in A}\delta_{n}^{2}
|\tilde{a}_{n+1}-\tilde{a}_{n}|^{2}+(1+\delta_{n}^{2})
|\tilde{a}_{n+1}+\tilde{a}_{n}|^{2}+\sum_{n\in
B}|\tilde{b}_{n}|^{2}\leq C \|U\|_{*}^{2}.
\end{align*}
On the other hand,
\begin{align*}
&\|U\|^{2}_{*}\leq 2\sum_{n\in
A}\delta_{n}^{2}\left(|\tilde{a}_{n}|^{2}
+|\tilde{a}_{n+1}|^{2}\right)
+|\tilde{a}_{n+1}+\tilde{a}_{n}|^{2}+\sum_{n\in
B}|\tilde{b}_{n}|^{2}\leq C'\|U\|^{2}_{\mathcal{Y}}.
\end{align*}
\end{Proof}\\
The following theorem provides a precise characterization of the
space $\mathcal{Y}$:
\begin{theo}\label{xx}
The space $\mathcal{Y}$ coincides algebraically and topologically
with the subspace of $\mathcal{H}_{0}\times \mathcal{H}_{-1}$
constituted by the
elements $((u^{0},v^{0},z^{0}),(u^{1},v^{1},z^{1}))$ such that\\
$(v^{0},v^{1})\in \mathcal{V}_{2}\times L^{2}(0,1)$, where
$\mathcal{V}_{2}$ is defined by \eqref{w3}.
\end{theo}
In order to prove this theorem, we need some results:
\begin{lemm}
For large $\lambda>0$, we have the following estimate
\begin{equation}\label{v1}
\begin{cases}
\tilde{v}_{\la}(x,\la)=
a_{2}(\rho_{2}(x)\sigma_{2}(x))^{-\frac{1}{4}}
\int_{x}^{1}\sqrt{\frac{\rho_{2}(t)}{\sigma_{2}(t)}}dt
\frac{\cos(\sqrt{\la}\int_{x}^{1}
\sqrt{\frac{\rho_{2}(t)}{\sigma_{2}(t)}}dt)}{2\la}
+\mathcal{O}(\frac{1}{\la^{\frac{3}{2}}}),\\
\frac{\partial\tilde{v}_{\la}}{\partial x}(x,\la)=a_{2}
(\rho_{2}(x))^{\frac{1}{4}}(\sigma_{2}(x))^{-\frac{3}{4}}
\int_{x}^{1}\sqrt{\frac{\rho_{2}(t)}{\sigma_{2}(t)}}dt
\frac{\sin(\sqrt{\la}\int_{x}^{1}
\sqrt{\frac{\rho_{2}(t)}{\sigma_{2}(t)}}dt)}{2\sqrt{\la}}
+\mathcal{O}(\frac{1}{\la}).
\end{cases}
\end{equation}
\end{lemm}
\begin{Proof}
By use of Liouville transformation (e.g., see \cite[Chapter 1]{BI})
and\\
$w(x)=\int_{x}^{1}\sqrt{\frac{\rho_{2}(t)}{\sigma_{2}(t)}}dt$
together with Lemma 2.1 in \cite{BI}, the solution
$\tilde{v}(x,\la)$ of Problem \eqref{w2} satisfies the following
integral equation
\begin{align}
&\tilde{v}(x,\la)=a_{2}\frac{(\rho_{2}(x)\sigma_{2}(x))^{-\frac{1}{4}}}
{\sqrt{\lambda}}
\sin(\sqrt{\lambda}w(x))\nonumber\\
&+\frac{(\rho_{2}(x)\sigma_{2}(x))^{-\frac{1}{4}}}{\sqrt{\lambda}}
\int_{0}^{w(x)}\sin(\sqrt{\lambda}(w(x)-t))Q(t)(\rho_{2}(t)
\sigma_{2}(t))^{\frac{1}{4}}\tilde{v}(t,\lambda)dt,\label{c2}
\end{align}
% $Q(w(x))=\frac{q_{2}(x)}{\rho_{2}(x)}
%-\left(\frac{\sigma_{2}(x)}{\rho_{2}^{3}(x)}\right)^{\frac{1}{4}}\left(\sigma_{2}(x)
%\left[(\sigma_{2}(x)\rho_{2}(x))^{-\frac{1}{4}}\right]'\right)'$.
where $Q(w)=
\frac{q_{2}(x)}{\rho_{2}(x)}-\frac{\ddot{\theta}(w)}{\theta(w)}$,
with $\theta(w)=(\rho_{2}(x)\sigma_{2}(x))^{\frac{1}{4}}$ and
$\dot{\theta}\equiv\frac{\partial \theta}{\partial w}$.
Differentiating this equation with respect to $\lambda$, we get
\begin{align}
&\tilde{v}_{\la}(x,\la)=a_{2}(\rho_{2}(x)\sigma_{2}(x))^{-\frac{1}{4}}
\left(w(x) \frac{\cos(\sqrt{\lambda}w(x))}{2\lambda}+
\frac{\sin(\sqrt{\lambda}w(x))}{2\lambda^{\frac{3}{2}}}\right)\nonumber\\
&+\frac{(\rho_{2}(x)\sigma_{2}(x))^{-\frac{1}{4}}} {2\lambda}
\int_{0}^{w(x)}
\sin(\sqrt{\lambda}(w(x)-t))Q(t)(\rho_{2}(t)\sigma_{2}(t))^{\frac{1}{4}}
\left(\tilde{v}_{\la}(t,\la)-\frac{\tilde{v}(t,\lambda)}{\sqrt{\lambda}}\right)dt\nonumber\\
%&-\frac{(\rho_{2}(x)\sigma_{2}(x))^{-\frac{1}{4}}}{\gamma_{2}\sqrt{\lambda}}\int_{x}^{1}\sin(\sqrt{\lambda}(w(x)-t))Q(t)v_{\la}(t,\lambda)dt\nonumber\\
&+\frac{(\rho_{2}(x)\sigma_{2}(x))^{-\frac{1}{4}}} {2\lambda
}\int_{0}^{w(x)}(w(x)-t)
\cos(\sqrt{\lambda}(w(x)-t))Q(t)(\rho_{2}(t)
\sigma_{2}(t))^{\frac{1}{4}}\tilde{v}(t,\lambda)dt.\label{c3}
\end{align}
%where $a(x,\lambda)=\frac{(\rho_{2}(x)\sigma_{2}(x))^{-\frac{1}{4}}}
%{2\lambda
%\gamma_{2}^{2}(\rho_{2}(1))^{\frac{-1}{4}}(\sigma_{2}(1))^{\frac{3}{4}}}$.
Let $M(\la)=\max\limits_{0\leq x \leq1}|\tilde{v}_{\la}(x,\la)|$.
Using \eqref{c2} and \eqref{c3}, one obtains for large $\lambda$,\\
$M(\la)\left(1-\frac{\int_{0}^{1}|Q(t)|dt}{\lambda}\right)\leq
\frac{C}{\lambda}$ $(c>0)$, and hence
\begin{equation}\label{c1}
|\tilde{v}_{\la}(x,\la)|=\mathcal{O}(\frac{1}{\lambda}).
\end{equation}
By substituting \eqref{g22} and \eqref{c1} into the integrals on the
right side of \eqref{c3}, the first estimate in \eqref{v1} holds.
The second estimate can be obtained in a similar
way.~~~~~~~~~~~~~~~~~~~~~~~~~~
\end{Proof}\\
We define the subsets
\begin{align}
\Lambda^{*}&=\left\{\hbox{$n$ and $-n$ such that}~n\in \Lambda\right\},\nonumber\\
&=\left\{\hbox{$n$ and $-n$ such
that}~\tilde{u}(0,\mu_{n})=\tilde{v}(0,\mu_{n})=0\right\},\label{h11}
\end{align}
\begin{equation}\label{b+}
B^{+}=\left\{n\in B~\hbox{such that $\mu_{n-1}$ is an eigenvalue of
Problem \eqref{g2}}~\right\}
\end{equation}
and
\begin{equation}\label{b-}
B^{-}=\left\{n\in B~\hbox{such that $\mu_{n-1}$ is an eigenvalue of
Problem \eqref{g1}}~\right\},
\end{equation}
where the set $\Lambda$ is defined in Proposition \ref{ff} and $B$
is defined by \eqref{aa2}. According to \eqref{g19}, we have
$B=B^{+}\cup B^{-}$.\\
%According to assertion $(i)$ of Lemma \ref{Lem1} and Corollary
%\ref{cor3} together with \eqref{aa1} and \eqref{h11}, we have the
%following result:
%\begin{lemm}\label{Lem2}
%If $n\in A\setminus \Lambda^{*}$, then $(\sqrt{\la_{n}})$ converges
%necessarily to an eigenvalue of Problem \eqref{g1} and an eigenvalue
%of Problem \eqref{g2}. Moreover, we have for some $C>0$
%$$\sqrt{\la_{n}}-\sqrt{\mu_{n-1}}\leq \frac{C}{n}~\hbox{and}~
%\sqrt{\mu_{n}}-\sqrt{\la_{n}}\leq \frac{C}{n},$$ with $\mu_{n}$ is
%an eigenvalue of Problem \eqref{g1} and $\mu_{n-1}$ is an eigenvalue
%of Problems \eqref{g2}, or conversely.
%\end{lemm}
We establish now the asymptotic behavior of the first derivative of
the eigenfunctions $(\phi_{n}(x))_{n\in \mathbb{Z}^{*}}$ at $x=1$.
\begin{prop}\label{prop1h}
For every $n\in\mathbb{Z}^*$, $\phi_{n}'(1)\neq0$. Furthermore, for
large $n$, \begin{equation}\label{am19}
 \left|\phi_{n}'(1)\right|=
\mathcal{O}(1),~n\in\Lambda^{*}\cup B^{+},
\end{equation}
\begin{equation}\label{am20}
\left|\phi_{n}'(1)\right|\leq\frac{C}{n},~C>0,~n\in A\setminus
\Lambda^{*},
\end{equation}
 \begin{equation}\left|\phi_{n+1}'(1)\right|=
\mathcal{O}(\frac{1}{n}),~ n\in A\label{am21}\end{equation}  and
\begin{equation} \left|\phi_{n}'(1)\right|=\mathcal{O}(\frac{1}{n}),
~n\in B^{-}.\label{am22}\end{equation}
\end{prop}
\begin{Proof}
It is clear from the expressions \eqref{ss} and \eqref{gg} of
${\phi}_{n}$ together with the initial conditions in \eqref{g7},
that
\begin{equation}\label{am10} \phi_{n}'(1)=\left\{
\begin{array}{lll}
-\sigma_{1}(0)\tilde{ u}_{x}(0,\la_{n}),& n\in\Lambda^*,\\
-{\sqrt{\la_n}\tilde{ u}(0,\la_n)}, &n \in
\mathbb{Z}^{*}\backslash\Lambda^*.
\end{array}
\right.
\end{equation}
Using this and \eqref{sig1}, the estimate \eqref{am19}
holds for $n\in \Lambda^{*}$.\\
In view of \eqref{v0}, \eqref{vv0} and \eqref{am10}, we obtain the
estimates
\eqref{am20} and \eqref{am21}.\\
Now by \eqref{la}, \eqref{g21} and \eqref{g19}, we have
\begin{align}
\sqrt{\lambda_{n}}\tilde{u}(0,\lambda_{n})\sim
\frac{\left(\rho_{1}(0)\sigma_{1}(0)\right)^{-\frac{1}{4}}}
{\left(\rho_{1}(-1)\right)^{\frac{1}{4}}
\left(\sigma_{1}(-1)\right)^{-\frac{3}{4}}},~\hbox{for large}~n\in
B^{+}\label{ttt}
\end{align}
and
\begin{equation}\label{ttt1}
\sqrt{\lambda_{n}}\tilde{u}(0,\lambda_{n})=\mathcal{O}(\frac{1}{n}),~\hbox{for
large}~n\in B^{-}.
\end{equation}
From this and \eqref{am10}, the estimates \eqref{am19} and
\eqref{am22} follow.
\end{Proof}\\
We are now in position to prove Theorem \ref{xx}. The method used in
the following proof was in some parts inspired
from the one of Theorem 5.4 in \cite{CAS}.\\
\begin{proof} {\bf  of Theorem \ref{xx}.}
To simplify notations, along the proof we denote the norm
$\|.\|_{\mathcal{V}_{2}\times L^{2}(0,1)}$ by $\|.\|_{H^{+}}$ and we
set $w(x)=\int_{x}^{1}\sqrt{\frac{\rho_{2}(t)}{\sigma_{2}(t)}}dt$.
First, we consider an element $U\in \mathcal{Y}$ and we prove that
$U_{\big|(0,1)}\in \mathcal{V}_{2}\times L^{2}(0,1)$. In view of
Proposition \ref{44}, any element $U\in \mathcal{Y}$ can be written
as
$$U=\sum_{n\in A}(c_{n}q_{n}
+d_{n}p_{n})+\sum_{n\in
B}\tilde{b}_{n}\frac{\overline{\Phi}_{n}}{\phi_{n}'(1)},$$ with the
coefficients $(\tilde{b}_{n})$, $(c_{n})$ and $(d_{n})$ are in
$\ell^{2}$. Clearly, we have
\begin{equation*}
\sum_{n\in A}c_{n}q_{n}+d_{n}p_{n}=\sum_{n\in A \setminus
\Lambda^{*}}c_{n}q_{n}+d_{n}p_{n}+\sum_{n\in
\Lambda^{*}}c_{n}q_{n}+d_{n}p_{n},
\end{equation*}
where the sets $A$ and $\Lambda^{*}$ are defined by \eqref{aa1} and
\eqref{h11}, respectively.\\
We set $U^{1,1}=\sum\limits_{n\in A\setminus
\Lambda^{*}}d_{n}p_{n}$. We prove that $\|U^{1,1}\|_{H^{+}}<\infty$.
From \eqref{gg} and \eqref{cc1}, one gets for $n\in A\setminus
\Lambda^{*}$
\begin{align}\label{l1}
{p_{n}}_{\big|(0,1)}=\frac{1}{2\delta_{n}}\left(
\frac{\sqrt{\la_{n+1}}\tilde{u}(0,\la_{n+1})}{\phi_{n+1}'(1)}
\overline{\tilde{v}}(x,\la_{n+1})-
\frac{\sqrt{\la_{n}}\tilde{u}(0,\la_{n})}{\phi_{n}'(1)}
\overline{\tilde{v}}(x,\la_{n})\right),
\end{align}
where $\overline{f}(x,\la_{n})=\left(f,i\sqrt{\lambda_{n}}f\right)$.
Substituting \eqref{am10} into \eqref{l1}, one has
\begin{align}\label{thh}
{p_{n}}_{\big|(0,1)}=-\frac{1}{2\delta_{n}}\left(
\overline{\tilde{v}}(x,\la_{n+1})-
\overline{\tilde{v}}(x,\la_{n})\right).
\end{align}
By use of the mean value theorem together with \eqref{g26} and
\eqref{v1}, we have
\begin{align}
\tilde{v}(x,\lambda_{n+1})-\tilde{v}(x,\lambda_{n})&
=(\sqrt{\la_{n+1}}-\sqrt{\la_{n}}) \frac{\partial}{\partial
\la}\tilde{v}(x,\la^{2})_{\big|\la=\sqrt{\alpha_{n}}},\nonumber\\
&=2\delta_{n}\sqrt{\alpha_{n}}\tilde{v}_{\lambda}
(x,\alpha_{n}),\label{tv}
\end{align}
for some $\sqrt{\alpha_{n}}\in (\sqrt{\la_{n}}, \sqrt{\la_{n+1}})$
and $n\in A$. Using this and \eqref{thh}, we obtain
\begin{align}\label{hj1}
{p_{n}}_{\big|(0,1)}=-\frac{1}{2\delta_{n}}
\left[2\delta_{n}\sqrt{\alpha_{n}}\tilde{v}_{\lambda}(x,\alpha_{n})
\left(1, i\sqrt{\la_{n}}\right) +\big(0,i\delta_{n}
\tilde{v}(x,\lambda_{n+1})\big)\right].
\end{align}
%On the other hand, it is easily seen from assertion $(i)$ of Lemma
%\ref{Lem1} and Corollary \ref{cor3}, that
%$(\sqrt{\lambda_{n}})_{n\in A\setminus \Lambda^{*}}$ converges
%necessarily to an eigenvalue of Problem \eqref{g2}, say
%$\sqrt{\mu_{n-1}}$. Clearly, we have
%$\sqrt{\alpha_{n}}=\sqrt{\mu_{n-1}}+\mathcal{O}(\frac{1}{n})$. From
%this and \eqref{v1}, we obtain
%\begin{align}
%&\sqrt{\alpha_{n}}\tilde{v}_{\la}(x,\alpha_{n})=\frac{
%a_{2}(\rho_{2}(x)\sigma_{2}(x))^{-\frac{1}{4}}
%w(x)}{2\sqrt{\alpha_{n}}}\left(\cos(\sqrt{\mu_{n-1}}w(x))
%+\mathcal{O}(\frac{1}{n})w(x)\sin(\sqrt{\mu_{n-1}}w(x))\right)\nonumber\\
%&+\mathcal{O}(\frac{1}{
%n^{2}})=\frac{\mu_{n-1}}{\sqrt{\alpha_{n}}}\tilde{v}_{\la}(x,\mu_{n-1})
%+\mathcal{O}(\frac{1}{n})w(x)\frac{\sqrt{\mu_{n-1}}}{\sqrt{\alpha_{n}}}
%\tilde{v}(x,\mu_{n-1}).\label{hj2}
%\end{align}
In view of \eqref{g22}, it is easily seen that $\left\|
\tilde{v}(x,\la_{n+1})\right\|^{2}_{L^{2}(0,1)}=\mathcal{O}(\frac{1}{n^{2}})$.
Hence
\begin{align}
&\left\|\sum_{n\in A\setminus \Lambda^{*}}\frac{d_{n}}{2\delta_{n}}
\delta_{n}\tilde{v}(x,\lambda_{n+1}) \right\|_{L^{2}(0,1)}^{2}\leq C
\left(\sum_{n\in A\setminus
\Lambda^{*}}|d_{n}|^{2}\right)\left(\sum_{n\in A\setminus
\Lambda^{*}}
\left\|\tilde{v}(x,\lambda_{n+1})\right\|_{L^{2}(0,1)}^{2}\right)\nonumber\\
&\leq C' \sum_{n\in A\setminus
\Lambda^{*}}|d_{n}|^{2}<\infty\label{bv}.
\end{align}
Using \eqref{v1}, it is clear that
$\left\|\sqrt{\alpha_{n}}\overline{\tilde{v}}_{\lambda}(x,\alpha_{n})
\right\|_{H^{+}}^{2}=\mathcal{O}(1)$. Thus, since
$\sqrt{\alpha_{n}}> \sqrt{\la_{n}}$, we have
\begin{align}
&\left\|\sum_{n\in A\setminus \Lambda^{*}}\frac{d_{n}}{\delta_{n}}
\sqrt{\alpha_{n}}\delta_{n}
\tilde{v}_{\lambda}(x,\alpha_{n})\left(1,i\sqrt{\la_{n}}\right)\right\|_{H^{+}}^{2}\leq
C \left\|\sum_{n\in A\setminus \Lambda^{*}}\frac{d_{n}}{\delta_{n}}
\sqrt{\alpha_{n}}\delta_{n}
\overline{\tilde{v}}_{\lambda}(x,\alpha_{n})\right\|_{H^{+}}^{2}\nonumber\\
&\leq C' \sum_{n\in A\setminus
\Lambda^{*}}|d_{n}|^{2}+C''\sum_{\underset{n\neq m}{ n,m\in
A\setminus \Lambda^{*}}}\left|d_{n}d_{m}\langle
\sqrt{\alpha_{n}}\overline{\tilde{v}}_{\la}(x,\alpha_{n}),
\sqrt{\alpha_{m}}\overline{\tilde{v}}_{\la}(x,\alpha_{m})\rangle_{H^{+}}\right|.
\label{bg1}
\end{align}
As noted before, for large $n\in A$, $\sqrt{\mu_{n-1}}$ close to
$\sqrt{\mu_{n}}$, and hence, by Corollary \ref{cor3},
\begin{equation}\label{yy1}
\sqrt{\la_{n}}-\sqrt{\mu_{n-1}}\leq \frac{C}{n}~\hbox{and}~
\sqrt{\mu_{n}}-\sqrt{\la_{n}}\leq \frac{C}{n},
\end{equation} for
some $C>0$. Since $\sqrt{\alpha_{n}}\in (\sqrt{\la_{n}},
\sqrt{\la_{n+1}})$ and
$\delta_{n}=\sqrt{\la_{n+1}}-\sqrt{\la_{n}}=\mathcal{O}(\frac{1}{n})$,
then
\begin{equation*}
|\sqrt{\alpha_{n}}-\sqrt{\mu_{n-1}}|\leq \frac{C}{n}~\hbox{and}~
|\sqrt{\alpha_{n}}-\sqrt{\mu_{n}}|\leq \frac{C}{n}.
\end{equation*}
Without loss of generality, we may assume that $\mu_{n-1}$ is an
eigenvalue of Problem \eqref{g2}. Then by \eqref{la}, $
|\sqrt{\alpha_{n}}-\frac{(n-1)\pi}{\gamma_{2}}|=\mathcal{O}
(\frac{1}{n})$. From this and \eqref{v1}, we have for large $n$,
\begin{equation*}
\sqrt{\alpha_{n}}\frac{\partial}{\partial
x}\tilde{v}_{\la}(x,\alpha_{n})=\frac{a_{2}}{2}
(\rho_{2}(x))^{\frac{1}{4}}(\sigma_{2}(x))^{-\frac{3}{4}}w(x)
\sin\left(\frac{(n-1)\pi}{\gamma_{2}}\right)
+\mathcal{O}(\frac{1}{n}).
\end{equation*}
Using this, yields
\begin{align}
&\left|\langle \sqrt{\alpha_{n}}\tilde{v}_{\la}(x,\alpha_{n}),
\sqrt{\alpha_{m}}\tilde{v}_{\la}(x,\alpha_{m})\rangle_{\mathcal{V}_{2}}
\right|=\left|\int_{0}^{1} \sqrt{\alpha_{n}}
\sqrt{\alpha_{m}}\frac{\partial}{\partial
x}(\tilde{v}_{\la})\frac{\partial}{\partial
x}((x,\alpha_{n})\tilde{v}_{\la})
(x,\alpha_{m})dx \right|\nonumber\\
&\leq
C\left|\int_{0}^{1}(\rho_{2}(x))^{\frac{1}{2}}(\sigma_{2}(x))^{
-\frac{3}{2}}
(w(x))^{2}\sin\left(\frac{(n-1)\pi}{\gamma_{2}}w(x)\right)
\sin\left(\frac{(m-1)\pi}{\gamma_{2}}w(x)\right)dx\right|\nonumber\\
&+\mathcal{O}\left(\frac{1}{|n|}\right)\left|\int_{0}^{1}(\rho_{2}(x))^{\frac{1}{4}}(\sigma_{2}(x))^{
-\frac{3}{4}} w(x)\sin\left(\frac{(m-1)\pi}{\gamma_{2}}w(x)\right)dx
\right|\nonumber\\
&+\mathcal{O}\left(\frac{1}{|m|}\right)\left|\int_{0}^{1}(\rho_{2}(x))^{\frac{1}{4}}
(\sigma_{2}(x))^{ -\frac{3}{4}}
w(x)\sin\left(\frac{(n-1)\pi}{\gamma_{2}}w(x)\right)dx
\right|\nonumber
\end{align}
Integrating by parts, we get
\begin{align}
&\mathcal{O}\left(\frac{1}{|n|}\right)\left|\int_{0}^{1}(\rho_{2}(x))^{\frac{1}{4}}(\sigma_{2}(x))^{
-\frac{3}{4}}w(x)\sin\left(\frac{(m-1)\pi}{\gamma_{2}}w(x)\right)dx
\right|\nonumber\\
&\leq
\mathcal{O}\left(\frac{1}{|n|}\right)\frac{(\rho_{2}(0)\sigma_{2}(0))
^{-\frac{1}{4}}\gamma_{2}^{2}}{(m-1)\pi}
+\mathcal{O}\left(\frac{1}{|nm|}\right)\leq \frac{C}{|nm|}.\nonumber
\end{align}
Similarly, twice integration by parts, yields
\begin{align*}
&\left|\int_{0}^{1}(\rho_{2}(x))^{\frac{1}{2}}(\sigma_{2}(x))^{
-\frac{3}{2}}
(w(x))^{2}\sin\left(\frac{(n-1)\pi}{\gamma_{2}}w(x)\right)
\sin\left(\frac{(m-1)\pi}{\gamma_{2}}w(x)\right)dx\right|\\
&\leq \frac{1}{2}\left|\int_{0}^{1}(\rho_{2}(x))^{\frac{1}{2}}
(\sigma_{2}(x))^{-\frac{3}{2}}
(w(x))^{2}\cos\left(\frac{(n-m)\pi}{\gamma_{2}}w(x)\right)dx\right|\nonumber\\
&+\frac{1}{2}\left|\int_{0}^{1}(\rho_{2}(x))^{\frac{1}{2}}
(\sigma_{2}(x))^{-\frac{3}{2}} (w(x))^{2}\cos\left(\frac{(n+m-2)\pi}
{\gamma_{2}}w(x)\right)dx\right|\nonumber\\
&\leq \frac{\gamma_{2}}{2|n-m|\pi}\left|\int_{0}^{1}
(\sigma_{2}^{-1}(x)
(w(x))^{2})'\sin\left(\frac{(n-m)\pi}{\gamma_{2}}w(x)\right)dx\right|\nonumber\\
&+ \frac{\gamma_{2}}{2|n+m-2|\pi}\left|\int_{0}^{1}
(\sigma_{2}^{-1}(x)
(w(x))^{2})'\sin\left(\frac{(n+m-2)\pi}{\gamma_{2}}w(x)\right)dx\right|\nonumber\\
%&\leq \left(
%\frac{\gamma_{2}^{2}}{2(n-m)^{2}\pi^{2}}
%+\frac{\gamma_{2}^{2}}{2(n+m-2)^{2}\pi^{2}}\right)\left(\frac{\sigma_{2}'(0)
%\gamma_{2}^{2}}
%{\rho_{2}(0)^{\frac{1}{2}}\sigma_{2}(0)^{\frac{3}{2}}}+2\gamma_{2}\sigma_{2}^{-1}(0)\right)\nonumber\\
%&+ \frac{\gamma_{2}^{2}}{2(n-m)^{2}\pi^{2}}\int_{0}^{1}\left|
%\left(\frac{(\sigma_{2}^{-1}(x) (w(x))^{2})'}{w'(x)}\right)'
%\cos\left(\frac{(n-m)\pi}{\gamma_{2}}w(x)\right)\right|dx\nonumber\\
%&+ \frac{\gamma_{2}^{2}}{2(n+m-2)^{2}\pi^{2}}\int_{0}^{1}\left|
%\left(\frac{(\sigma_{2}^{-1}(x) (w(x))^{2})'}{w'(x)}\right)'
%\cos\left(\frac{(n-m)\pi}{\gamma_{2}}w(x)\right)\right|dx\nonumber\\
&\leq \frac{C'}{(m-n)^{2}} +\frac{C''}{(m+n)^{2}}\leq
\frac{C}{(m-n)^{2}}.
\end{align*}
%Analogously, we have
%\begin{align*}
%&\left|\langle \alpha_{n}\tilde{v}_{\la}(x,\alpha_{n}),
%\alpha_{m}\tilde{v}_{\la}(x,\alpha_{m})\rangle_{L^{2}(0,1)}
%\right|=\left|\int_{0}^{1}\alpha_{n}
%\alpha_{m}\tilde{v}_{\la}(x,\alpha_{n})\tilde{v}_{\la}
%(x,\alpha_{m})dx \right|\nonumber\\
%&\leq \frac{C_{1}}{(n-m)^{2}}+\frac{C_{2}}{(n+m)^{2}}.
%\end{align*}
%Then from this, \eqref{fk} and \eqref{fk1}, to prove that
%\begin{equation}\label{fkk}
%\sum_{\underset{n\neq m}{ n,m\in A\setminus
%\Lambda^{*}}}\left|d_{n}d_{m}\langle
%\sqrt{\alpha_{n}}\tilde{v}_{\la}(x,\alpha_{n}),
%\sqrt{\alpha_{m}}\tilde{v}_{\la}(x,\alpha_{m})\rangle_{\mathcal{V}_{2}}
%\right|< \infty,
%\end{equation}
Analogously, we have $\left|\langle
\alpha_{n}\tilde{v}_{\la}(x,\alpha_{n}),
\alpha_{m}\tilde{v}_{\la}(x,\alpha_{m})\rangle_{L^{2}(0,1)}
\right|\leq \frac{C'}{(n-m)^{2}}.$
%\begin{align*}
%\left|\langle \alpha_{n}\tilde{v}_{\la}(x,\alpha_{n}),
%\alpha_{m}\tilde{v}_{\la}(x,\alpha_{m})\rangle_{L^{2}(0,1)}
%\right|\leq \frac{C'}{(n-m)^{2}}.
%\end{align*}
Hence from the above, in order to prove that the last term in
\eqref{bg1} is finite, it suffices to show that
\begin{equation*}
\sum_{\underset{n\neq m}{ n,m\in A\setminus
\Lambda^{*}}}\frac{\left|d_{n}d_{m}\right|}{(n-m)^{2}}< \infty.
\end{equation*}
Clearly,
\begin{align*}
&\sum_{\underset{n\neq m}{ n,m\in A\setminus
\Lambda^{*}}}\frac{\left|d_{n}d_{m}\right|}{(n-m)^{2}} \leq
\sum_{j\neq 0}\frac{1}{j^{2}}\sum_{n\in A\setminus
\Lambda^{*}}\left|d_{n}d_{n+j}\right|\nonumber\\
&\leq \sum_{j\neq 0}\frac{1}{j^{2}}\sum_{n\in A\setminus
\Lambda^{*}}\left|d_{n}\right|^{2}\leq C \sum_{n\in A\setminus
\Lambda^{*}}\left|d_{n}\right|^{2}<\infty.
\end{align*}
Therefore,
\begin{align}
&\sum_{\underset{n\neq m}{ n,m\in A\setminus
\Lambda^{*}}}\left|d_{n}d_{m}\langle
\sqrt{\alpha_{n}}\overline{\tilde{v}}_{\la}(x,\alpha_{n}),
\sqrt{\alpha_{m}}\overline{\tilde{v}}_{\la}(x,\alpha_{m})\rangle_{H^{+}}
\right|\leq C \sum_{ n\in A\setminus \Lambda^{*}}
\left|d_{n}\right|^{2}<\infty.\label{bg2}
\end{align}
Combining this with \eqref{bv} and \eqref{bg1}, we conclude that
\begin{equation}\label{rv}
\|U^{1,1}\|_{H^{+}}^{2}=\left\|\sum_{n\in A\setminus
\Lambda^{*}}d_{n}p_{n}\right\|_{H^{+}}^{2}\leq C \sum_{n\in A
\setminus \Lambda^{*}}|d_{n}|^{2}.
\end{equation}
Let $U^{1,2}=\sum\limits_{n\in A\setminus \Lambda^{*}} c_{n}q_{n}$.
Then
\begin{equation*}
\|U^{1,2}\|_{H^{+}}^{2} \leq 2\left\|\sum_{n \in A\setminus
\Lambda^{*}}c_{n}\left(q_{n}-\frac{\overline{\Phi}_{n}}{\phi_{n}'(1)}\right)
\right\|_{H^{+}}^{2}+2\left\|\sum_{n \in A\setminus
\Lambda^{*}}c_{n}\frac{\overline{\Phi}_{n}}{\phi_{n}'(1)}\right\|_{H^{+}}^{2}.
\end{equation*}
In view of \eqref{cc1} and \eqref{rv}, one has
\begin{equation}\label{th}
\left\|\sum_{n\in A\setminus
\Lambda^{*}}c_{n}\left(q_{n}-\frac{\overline{\Phi}_{n}}{\phi_{n}'(1)}
\right)\right\|_{H^{+}}^{2}=\left\|\sum_{n\in A\setminus
\Lambda^{*}}c_{n}\delta_{n}p_{n}\right\|_{H^{+}}^{2}\leq
C'\sum_{n\in A\setminus \Lambda^{*}}|c_{n}|^{2}<\infty.
\end{equation}
By \eqref{gg} and \eqref{am10},
\begin{equation}\label{yy}
\frac{\overline{\Phi}_{n}}{\phi_{n}'(1)}=\overline{\tilde{v}}(x,\lambda_{n}),~~~
\hbox{for}~n\in A \setminus \Lambda^{*}.
\end{equation}
%By Corollary \ref{cor3}, we have for large $n\in A\setminus
%\Lambda^{*}$
%\begin{equation}\label{yy1}
%\sqrt{\la_{n}}=\sqrt{\mu_{n-1}}+\xi_{n},~\hbox{ with}~ \xi_{n}\leq
%\frac{C}{n}.
%\end{equation}
Then analogously to \eqref{tv}, using \eqref{yy1}, we get for some
$\sqrt{\beta_{n}}\in (\sqrt{\mu_{n-1}}, \sqrt{\la_{n}})$,
\begin{align}
\tilde{v}(x,\la_{n})-\tilde{v}(x,\mu_{n-1})=2\xi_{n}\sqrt{\beta_{n}}
\tilde{v}_{\lambda}(x,\beta_{n}),\label{tv1}
\end{align}
with $|\xi_{n}|\leq \frac{C}{n}$ for large $n$. Substituting
\eqref{tv1} into \eqref{yy}, yields
\begin{equation}\label{tv2}
\frac{\overline{\Phi}_{n}}{\phi_{n}'(1)}=
\overline{\tilde{v}}(x,\mu_{n-1})+2\xi_{n}\sqrt{\beta_{n}}
\tilde{v}_{\lambda}
(x,\beta_{n})(1,i\sqrt{\mu_{n-1}})+\left(0,\xi_{n}\tilde{v}(x,\la_{n})\right).
\end{equation}
By \eqref{v1}, $\|\tilde{v}_{\lambda}
(x,\beta_{n})(1,i\sqrt{\mu_{n-1}})\|^{2}_{H^{+}}=\mathcal{O}(\frac{1}{n^{2}})$.
Furthermore, under the above assumption about $\mu_{n-1}$, since
$\{\overline{\tilde{v}}(x,\mu_{n})~\hbox{such
that}~\tilde{v}(0,\mu_{n})= 0\}_{n\in \mathbb{Z}^{*}}$ is an
orthogonal subset in $H^{+}$ with
$\|\overline{\tilde{v}}(x,\mu_{n})\|^{2}_{H^{+}}=\mathcal{O}(1)$,
then from \eqref{g22} and \eqref{tv2}, we have
\begin{align}
&\left\|\sum_{n\in A\setminus
\Lambda^{*}}c_{n}\frac{\overline{\Phi}_{n}}{\phi_{n}'(1)}\right\|_{H^{+}}^{2}\leq
C_{1}\left\|\sum_{n\in A\setminus
\Lambda^{*}}c_{n}\xi_{n}\tilde{v}(x,\la_{n})\right\|_{L^{2}(0,1)}^{2}
+C_{2} \left\|\sum_{n\in A\setminus
\Lambda^{*}}c_{n}\overline{\tilde{v}}(x,\mu_{n-1})\right\|_{H^{+}}^{2}\nonumber\\
&+C_{3}\left\|\sum_{n\in A\setminus
\Lambda^{*}}c_{n}\tilde{v}_{\la}(x,\beta_{n})(1,i\sqrt{\mu_{n-1}})
\right\|_{H^{+}}^{2}\leq C_{1}\left(\sum_{n\in A\setminus
\Lambda^*}|c_{n}\xi_{n}|^{2}\right)\sum_{n\in A\setminus
\Lambda^*}\|\tilde{v}(x,\la_n)\|^{2}_{L^{2}(0,1)}\nonumber\\
&+C'_{2}\sum_{n\in A\setminus
\Lambda^*}|c_{n}|^{2}+C_{3}\left(\sum_{n\in A\setminus
\Lambda^*}|c_{n}|^{2}\right)\sum_{n\in A\setminus
\Lambda^*}\|\tilde{v}_{\la}(x,\beta_n)(1,i\sqrt{\mu_{n-1}})\|^{2}_{H^+}
\nonumber\\
&\leq C' \sum_{n\in A\setminus \Lambda^{*}}|c_{n}|^{2}\label{s2}.
\end{align}
Hence, from this and \eqref{th}, it follows that
$\|U^{1,2}\|_{H^{+}}^{2}\leq C \sum\limits_{n\in A\setminus
\Lambda^{*}}|c_{n}|^{2}$. Therefore, $U^{1,1}_{\big|(0,1)}$ and
$U^{1,2}_{\big|(0,1)}$ belong to
$\mathcal{V}_{2}\times L^{2}(0,1)$.\\
Let us now take $U^{2,1}=\sum\limits_{n\in \Lambda^{*}}d_{n}p_{n}$.
Recall that if $n\in \Lambda^{*}$, then $\lambda_{n}=\mu_{n}$ and
$(n+1)\not\in \Lambda^{*}$, so that for $n\in \Lambda^{*}$,
$\Phi_{n}$ and $\Phi_{n+1}$ have the forms \eqref{ss} and
\eqref{gg}, respectively. Proceeding as above, in view of
\eqref{cc1}, \eqref{am10} and \eqref{tv}, one gets
\begin{align*}
&p_{n}=\frac{1}{2\delta_{n}}\left(\frac{\sqrt{\la_{n+1}}
\tilde{u}(0,\la_{n+1})}{\phi_{n+1}'(1)}\overline{\tilde{v}}(x,\la_{n+1})-
\frac{\sigma_{1}(0) \tilde{u}_{x}(0,\la_{n})}{\phi_{n}'(1)}
\overline{\tilde{v}}(x,\la_{n})\right)\nonumber\\
&=-\frac{1}{2\delta_{n}}\left(\overline{\tilde{v}}
(x,\lambda_{n+1})-\overline{\tilde{v}}
(x,\lambda_{n})\right)\nonumber\\
&=-\frac{1}{2\delta_{n}}
\left(2\sqrt{\alpha_{n}}\delta_{n}\tilde{v}_{\lambda}
(x,\alpha_{n})(1,i\sqrt{\la_{n}})+
(0,i\delta_{n}\tilde{v}(x,\la_{n+1}))\right),
\end{align*}
for some $\sqrt{\alpha_n}\in (\sqrt{\la_{n}},\sqrt{\la_{n+1}})$. In
the same way as in \eqref{bv}, \eqref{bg1} and \eqref{bg2}, it can
be shown that
%From \eqref{ghh}, obviously $d_{n}$ can be written as
%$\delta_{n}\tilde{d}_{n}$ with $(\delta_{n}\tilde{d}_{n})\in
%\ell^{2}$ for $n\in \Lambda^{*}$. Hence, similarly to \eqref{bg1}
%and \eqref{bg2}, using \eqref{g22}, \eqref{v1} and \eqref{dr}, one
%gets
%\begin{align*}
%&\left\|\sum_{n\in \Lambda^{*}}\frac{d_{n}}{\delta_{n}}
%\sqrt{\alpha_{n}}\delta_{n}\tilde{v}_{\lambda}
%(x,\alpha_{n})(1,i\sqrt{\la_{n}})\right\|_{H^{+}}^{2}\leq
%\left\|\sum_{n\in
%\Lambda^{*}}d_{n}\sqrt{\alpha_{n}}\overline{\tilde{v}}_{\lambda}(
%x,\alpha_{n})\right\|_{H^{+}}^{2}\leq C' \sum_{n\in
%\Lambda^{*}}|d_{n}|^{2}.
%\end{align*}
%Using this, \eqref{g22} and \eqref{dr}, we obtain
\begin{align*}
&\|U^{2,1}\|_{H^{+}}^{2}=\left\|\sum_{n\in
\Lambda^{*}}d_{n}p_{n}\right\|_{H^{+}}^{2}\leq
%C_{1}\left\|\sum_{n\in
%\Lambda^{*}}d_{n}\sqrt{\alpha_{n}}\overline{\tilde{v}}_{\lambda}(
%x,\alpha_{n})\right\|_{H^{+}}^{2}\nonumber\\
%&+C_{2}\left\|\sum_{n\in \Lambda^{*}}d_{n}\tilde{v}(x,\lambda_{n+1})
%\right\|_{L^{2}(0,1)}^{2}\leq
C \sum_{n\in \Lambda^{*}}|d_{n}|^{2}.
\end{align*}
We set $U^{2,2}=\sum\limits_{n\in \Lambda^{*}} c_{n}q_{n}$.
Similarly, we have
\begin{align*}
&\|U^{2,2}\|_{H^{+}}^{2} \leq 2\left\|\sum_{n \in
\Lambda^{*}}c_{n}\left(q_{n}-\frac{\overline{\Phi}_{n}}{\phi_{n}'(1)}\right)
\right\|_{H^{+}}^{2}+2\left\|\sum_{n \in
\Lambda^{*}}c_{n}\frac{\overline{\Phi}_{n}}{\phi_{n}'(1)}
\right\|_{H^{+}}^{2}\nonumber\\
&\leq 2\left\|\sum_{n\in
\Lambda^{*}}c_{n}\delta_{n}p_{n}\right\|_{H^{+}}^{2}+2\left\|\sum_{n\in
\Lambda^{*}}c_{n}\overline{\tilde{v}}(x,\la_{n})\right
\|_{H^{+}}^{2}\leq C \sum_{n\in \Lambda^{*}}|c_{n}|^{2}.
\end{align*}
%Taking into account that
%$\{\overline{\tilde{v}}(x,\lambda_{n})\}_{n\in \Lambda^{*}}$ is an
%orthogonal subset in $H^{+}$ with\\
%$\|\overline{\tilde{v}}(x,\lambda_{n})\|^{2}_{H^{+}}=\mathcal{O}(1)$,
%then from \eqref{am10}, we have
%\begin{align}
%&\left\|\sum_{n\in
%\Lambda^{*}}c_{n}\frac{\overline{\Phi}_{n}}{\phi_{n}'(1)}\right\|_{H^{+}}^{2}\leq
%C_{1} \sum_{n\in  \Lambda^{*}}|c_{n}|^{2}\label{s22}.
%\end{align}
%Whence, by \eqref{th1} and \eqref{s22}, one gets
%\begin{equation*}
%\|U^{2,2}\|_{H^{+}}^{2} \leq C\sum_{n \in \Lambda^{*}}|c_{n}|^{2}.
%\end{equation*}
%In view of \eqref{ss} and taking into account that
%$\{\overline{\tilde{v}}(x,\lambda_{n})\}_{n\in \Lambda^{*}}$ is an
%orthogonal subset in $H^{1}(0,1)\times L^{2}(0,1)$, then from
%\eqref{y2} and \eqref{gb}, one has
%\begin{equation*}
%\|U^{2,2}\|_{H^{+}}=C\left(\sum_{n \in
%\Lambda^{*}}|c_{n}|^{2}\right)^{\frac{1}{2}}.
%\end{equation*}
Therefore, $U^{2,1}_{\big|(0,1)}$ and $U^{2,2}_{\big|(0,1)}$ belong
to $\mathcal{V}_{2}\times L^{2}(0,1)$.\\
Recall that $B=B^{+}\cup B^{-}$ and set $U^{3,1}=\sum\limits_{n\in
B^{+}}\tilde{b}_{n}\frac{\overline{\Phi}_{n}}{\phi_{n}'(1)}$.
Analogously to \eqref{tv2}, using \eqref{gg}, \eqref{g19} and
\eqref{am10}, one has
\begin{align*}
&\frac{\overline{\Phi}_{n}}{\phi_{n}'(1)}=\overline{\tilde{v}}(x,\mu_{n-1})+
2\mathcal{O}(\frac{1}{n})\sqrt{\beta_{n}} \tilde{v}_{\lambda}
(x,\beta_{n})(1,i\sqrt{\mu_{n-1}})+
\left(0,\mathcal{O}(\frac{1}{n})\tilde{v}(x,\la_{n})\right),
\end{align*}
for some $\sqrt{\beta_{n}}\in (\sqrt{\mu_{n-1}},\sqrt{\la_{n}})$.
Now, similarly to \eqref{s2}, taking into account that
$\{\overline{\tilde{v}}(x,\mu_{n-1})\}_{n\in B^{+}}$ is an
orthogonal subset in $H^{+}$ with
$\|\overline{\tilde{v}}(x,\mu_{n-1})\|^{2}_{H^{+}}=\mathcal{O}(1)$,
we obtain $\|U^{3,1}\|_{H^{+}}^{2}\leq C\sum\limits_{n\in B^{+}}
|\tilde{b}_{n}|^{2}$. Thus $U^{3,1}_{\big|(0,1)}\in \mathcal{V}_{2}\times L^{2}(0,1)$.\\
Let us define $U^{3,2}=\sum\limits_{n\in
B^{-}}\tilde{b}_{n}\frac{\overline{\Phi}_{n}}{\phi_{n}'(1)}$. Recall
that for $n\in B^{-}$, $\mu_{n-1}$ is not an eigenvalue of Problem
\eqref{g2}, i.e., $\tilde{v}(0,\mu_{n-1})\neq0$ and so
$\{\overline{\tilde{v}}(x,\mu_{n-1})\}_{n\in B^{-}}$ is not
necessarily an orthogonal subset in $H^+$. By \eqref{g22},
\eqref{gg} and \eqref{am10}, we have
$\|\frac{\overline{\Phi}_{n}}{\phi_{n}'(1)}\|_{H^{+}}=\mathcal{O}(1)$
for $n\in B^{-}$. Hence
\begin{align}\label{amm2}
&\|U^{3,2}\|_{H^{+}}^{2}\leq C\sum_{n\in B^{-}}
|\tilde{b}_{n}|^{2}+\sum_{\underset{n\neq m}{ n,m\in
B^{-}}}\left|\tilde{b}_{n}\tilde{b}_{m}\langle
\frac{\overline{\Phi}_{n}}{\phi_{n}'(1)},
\frac{\overline{\Phi}_{m}}{\phi_{m}'(1)}\rangle_{H^{+}}\right|.
\end{align}
%On the other hand, from \eqref{ghh}, \eqref{am10} and \eqref{ttt1},
From \eqref{ghh} and \eqref{am22}, it is easily seen that
$\tilde{b}_n\sim\frac{b_{n}}{n}$. On the other hand, clearly by
\eqref{pro2}, \eqref{xx1} and \eqref{xx2}, a simple computation
gives
\begin{equation}\label{vb}
b_n=\mathcal{O}(1),~\hbox{for large}~ n\in B^{-}.
\end{equation}
Since $(\tilde{b}_{n})\in \ell^{2}$, then to see that
$\sum\limits_{\underset{n\neq m}{ n,m\in
B^{-}}}\left|\tilde{b}_{n}\tilde{b}_{m}\langle
\frac{\overline{\Phi}_{n}}{\phi_{n}'(1)},
\frac{\overline{\Phi}_{m}}{\phi_{m}'(1)}\rangle_{H^{+}}\right|<\infty$,
it suffices to show that
\begin{equation}\label{amm1}
\sum_{\underset{n\neq m}{ n,m\in B^{-}}}\left|\frac{1}{nm}\langle
\frac{\overline{\Phi}_{n}}{\phi_{n}'(1)},
\frac{\overline{\Phi}_{m}}{\phi_{m}'(1)}\rangle_{H^{+}}\right|=
\sum_{\underset{n\neq m}{ n,m\in B^{-}}}\left|\frac{1}{nm}\langle
\overline{\tilde{v}}(x,\la_{n}),
\overline{\tilde{v}}(x,\la_{m})\rangle_{H^{+}}\right|<\infty.
\end{equation}
%By \eqref{gg} and \eqref{am10},
%\begin{equation}\label{j8}
%\frac{\overline{\Phi}_{n}}{\phi_{n}'(1)}=\overline{\tilde{v}}(x,\la_{n}).
%\end{equation}
Hence, using \eqref{g22} and integrating by parts, we get
\begin{align*}
&\left|\langle \overline{\tilde{v}}(x,\la_{n}),
\overline{\tilde{v}}(x,\la_{m})\rangle_{H^{+}}\right|
\leq\frac{C_{1}}{|n-m|}+\frac{C_{2}}{|n+m|}\leq\frac{C}{|n-m|}.
\end{align*}
Now, we show that $\sum\limits_{\underset{n\neq m}{n,m\in
\mathbb{Z}^{*}}}\frac{1}{|nm(n-m)|}< \infty$. For fixed $n$,
$\sum\limits_{\underset{n\neq m}{m\in
\mathbb{Z}^{*}}}\frac{1}{|mn(n-m)|}$ converges. Indeed for large
enough $N$ and fixed $n$, we have
\begin{equation*}
\sum_{\underset{n\neq
m}{m=1}}^{N}\frac{1}{mn(m-n)}=\sum_{m=1}^{n-1}\frac{1}{mn(n-m)}+\sum_{m=n+1}^{N}\frac{1}{mn(m-n)}.
\end{equation*}
It is clear that
\begin{align*}
&\sum_{m=n+1}^{N}\frac{1}{mn(m-n)}=\frac{1}{n^{2}}\left(\sum_{m=n+1}^{N}\frac{1}{m-n}-\sum_{m=n+1}^{N}\frac{1}{m}\right),\\
&=\frac{1}{n^{2}}\left(\sum_{j=1}^{N-n}\frac{1}{j}-
\sum_{j=n+1}^{N}\frac{1}{j}\right)=\frac{1}{n^{2}}\left(\sum_{j=1}^{n}\frac{1}{j}-\sum_{j=N-n+1}^{N}\frac{1}{j}\right)
\end{align*}
and $\lim\limits_{N\rightarrow
+\infty}\frac{1}{n^{2}}\left(\sum\limits_{j=1}^{n}\frac{1}{j}-\sum\limits_{j=N-n+1}^{N}
\frac{1}{j}\right)=\frac{1}{n^{2}}\sum\limits_{j=1}^{n}\frac{1}{j}$.
In a similar way, we have
\begin{align*}
&\sum_{m=1}^{n-1}\frac{1}{mn(n-m)}=\frac{1}{n^{2}}\left(
\sum_{m=1}^{n-1}\frac{1}{n-m}+\sum_{m=1}^{n-1}\frac{1}{m}\right)\\
&=\frac{1}{n^{2}}\left(
\sum_{j=1}^{n-1}\frac{1}{j}+\sum_{m=1}^{n-1}\frac{1}{m}\right)=\frac{2}{n^{2}}\sum_{j=1}^{n-1}\frac{1}{j}.
\end{align*}
Since $\sum\limits_{n\in \mathbb{Z}^{*}}\frac{1}{n^{2}}\left(
\sum\limits_{j=1}^{n}\frac{1}{j}\right)<\infty$, then the estimate
\eqref{amm1} follows. From this
and \eqref{amm2}, we have $\|U^{3,2}\|_{H^{+}}^{2}< \infty$. Therefore, $U^{3,2}_{\big|(0,1)} \in \mathcal{V}_{2}\times L^{2}(0,1)$.\\
We have proved that if $U\in \mathcal{Y}$, then the restriction of
$U$ to $(0,1)$ belongs to\\$\mathcal{V}_{2}\times L^{2}(0,1)$.\\
Consider now an element $U=((u^{0},v^{0},z^{0}),(u^{1},v^{1},z^{1}))
\in X_{-\frac{1}{2}}\times X_{-1}$ such that $U_{\big |(0,1)}\in
\mathcal{V}_{2}\times L^{2}(0,1)$. We shall prove that
$U\in \mathcal{Y}$.\\
%Since $U\in X_{-\frac{1}{2}}\times X_{-1}$ and in view of
%\eqref{hn}, it can be written as follows
Following the proof of Proposition \ref{44}, $U \in
X_{-\frac{1}{2}}\times X_{-1}$ can be written as follows
\begin{align}
&U=\sum_{n\in \mathbb{Z}^{*}}a_{n}\overline{\Phi}_{n}=\sum_{n\in
A}c_{n}q_{n}+d_{n}\delta_{n}p_{n}+\sum_{n\in
B}\tilde{b}_{n}\frac{\overline{\Phi}_{n}}{\phi_{n}'(1)},\label{jjj}
\end{align}
where $c_{n}=\tilde{a}_{n+1}+\tilde{a}_{n}$ and
$d_{n}=\tilde{a}_{n+1}-\tilde{a}_{n}$ for $n\in A$. According to
Proposition \ref{44}, to see that $U\in \mathcal{Y}$, it is
sufficient to show that $(c_{n})_{n\in A}\in \ell^{2}$,
$(\delta_{n}d_{n})_{n\in A}\in \ell^{2}$ and $(\tilde{b}_{n})_{n\in
B}\in \ell^{2}$.\\
%Since $U\in X_{-\frac{1}{2}}\times X_{-1}$, then from \eqref{xa},
%\eqref{k1} and \eqref{jjj}, we have
%\begin{align}
%&\sum_{n\in B}\frac{|b_{n}|^{2}}{\la_{n}}<\infty.\label{t1}
%\end{align}
Since $\tilde{b}_{n}=b_{n}\phi_{n}'(1)$, then by \eqref{am22} and
\eqref{vb}, we get $\sum\limits_{n\in
B^{-}}|\tilde{b}_{n}|^{2}<\infty$. Proceeding as above, since
$(\tilde{b}_{n})_{n\in B^{-}}\in \ell^{2}$, then
$U^{3,2}_{\big|(0,1)}=\sum\limits_{n\in
B^{-}}\tilde{b}_{n}\frac{\overline{\Phi}_{n}}{\phi_{n}'(1)}_{\big|(0,1)}\in
\mathcal{V}_{2}\times L^{2}(0,1)$. Therefore, $U^{*}= U-U^{3,2}$ is
also such that, when restricted to $(0,1)$, belongs to
$\mathcal{V}_{2}\times L^{2}(0,1)$. As above, using \eqref{ss} and
\eqref{gg}, one has
\begin{align}
&U^{*}_{\big|(0,1)}=\left(\sum\limits_{n\in \mathbb{Z}^{*}\setminus
B}a_{n}\overline{\Phi}_{n}+\sum\limits_{n\in
B^{+}}b_{n}\overline{\Phi}_{n}\right)_{\big|(0,1)}=-\sum_{n\in
\Lambda^{*}}\tilde{a}_{n}
\overline{\tilde{v}}(x,\lambda_{n})\nonumber\\
&-\sum_{n\in A \setminus \Lambda^{*} }\tilde{a}_{n}\left[
\left(\overline{\tilde{v}}(x,\mu_{n-1})+ 2\xi_{n}\sqrt{\alpha_{n}}
\tilde{v}_{\lambda}(x,\alpha_{n})(1,i\sqrt{\mu_{n-1}})\right)+
\xi_{n}\left(0,i\tilde{v}(x,\lambda_{n})\right)\right]\nonumber\\
&-\sum_{n\in A}\left[\tilde{a}_{n+1}
\left(\overline{\tilde{v}}(x,\mu_{n-1})+
2\mathcal{O}(\frac{1}{n})\sqrt{\alpha'_{n}}
\tilde{v}_{\lambda}(x,\alpha'_{n})(1,i\sqrt{\mu_{n-1}})\right)+
\mathcal{O}(\frac{1}{n})\left(0,i\tilde{v}(x,\lambda_{n+1})\right)\right]\nonumber\\
&-\sum_{n\in B^{+}}\tilde{b}_{n}\left[\left(
\overline{\tilde{v}}(x,\mu_{n-1})+
2\mathcal{O}(\frac{1}{n})\sqrt{\beta_{n}}
\tilde{v}_{\lambda}(x,\beta_{n})(1,i\sqrt{\mu_{n-1}})\right)+\mathcal{O}(\frac
{1}{n})\left(0,\tilde{v}(x,\lambda_{n})\right)\right],\label{f1}
\end{align}
for some $\sqrt{\alpha_{n}}\in (\sqrt{\mu_{n-1}},\sqrt{\la_{n}})$
for $n\in A \setminus \Lambda^{*}$, $\sqrt{\alpha'_{n}}\in
(\sqrt{\mu_{n-1}},\sqrt{\la_{n+1}})$ for $n\in A$ and
$\sqrt{\beta_{n}}\in (\sqrt{\mu_{n-1}},\sqrt{\la_{n}})$ for $n\in
B^{+}$. Since $U\in X_{-\frac{1}{2}}\times X_{-1}$, then from
\eqref{xa}, \eqref{k2} and \eqref{jjj}, we have
\begin{equation}\label{i22}
\sum_{n\in A\setminus \Lambda^{*}}\frac{|a_{n}|^{2}}{\la_{n}}\langle
\Phi_{n}, \Phi_{n} \rangle_{\mathcal{W}}<\infty.
\end{equation}
From assertion $(ii)$ of Lemma \ref{Lem1} and the estimates
\eqref{k5} and \eqref{k44} together with \eqref{am20} and
\eqref{am10}, we obtain $$\langle \Phi_{n}, \Phi_{n}
\rangle_{\mathcal{W}}\sim |\phi_{n}'(1)|^{2},~\hbox{for}~n\in
A\setminus \Lambda^{*}.$$ As $\tilde{a}_{n}=a_{n}\phi_{n}'(1)$, then
by \eqref{i22}, $\sum\limits_{n\in A\setminus
\Lambda^{*}}\frac{|\tilde{a}_{n}|^{2}}{n^{2}}<\infty$. Since
$|\xi_{n}|\leq \frac{C}{n}$ for large $n$, then
\begin{equation}\label{bg}
\sum_{n\in A\setminus \Lambda^{*}}|\tilde{a}_{n}\xi_{n}|^{2}<\infty.
\end{equation}
 Analogously, using \eqref{k1} and \eqref{k3} together with
\eqref{am19} and \eqref{am21}, we get
\begin{equation}\label{i2}
\sum_{n\in
\Lambda^{*}}\frac{|\tilde{a}_{n}|^{2}}{n^{2}}<\infty,~\sum_{n\in
B^{+}}\frac{|\tilde{b}_{n}|^{2}}{n^{2}}<\infty~
\hbox{and}~\sum_{n\in A}\frac{|\tilde{a}_{n+1}|^{2}}{n^{2}}< \infty.
\end{equation}
%In view of \eqref{am10} and \eqref{yy1} together with \eqref{i222}
%and \eqref{i2}, it is easily seen that
%\begin{equation}\label{jd}
%\sum\limits_{n\in A\setminus \Lambda^{*}}
%|\tilde{a}_{n}\xi_{n}|^{2}<\infty~\hbox{and}~\sum\limits_{n\in A}
%|\tilde{a}_{n+1} \mathcal{O}(\frac{1}{n})|^{2}<\infty.
%\end{equation}
As $\delta_{n}=\mathcal{O}(\frac{1}{n})$ for $n\in A$, then from
\eqref{bg} and \eqref{i2}, we conclude that
$(\delta_{n}d_{n})_{n\in A}\in \ell^{2}$.\\
In view of \eqref{g22}, \eqref{bg} and \eqref{i2},
\begin{equation*}\label{f2}
\left\|\sum\limits_{n\in A\setminus \Lambda^{*}}\tilde{a}_{n}\xi_{n}
\tilde{v}(x,\lambda_{n})\right\|_{L^{2}(0,1)}\leq C\sum\limits_{n\in
A\setminus \Lambda^{*}} \left|\tilde{a}_{n}\xi_{n}\right|^{2}<
\infty,
\end{equation*}
\begin{equation*}\label{f22}
\left\|\sum\limits_{n\in A}\tilde{a}_{n+1} \mathcal{O}(\frac{1}{n})
\tilde{v}(x,\lambda_{n+1})\right\|_{L^{2}(0,1)} \leq C
\sum\limits_{n\in A}
\left|\frac{\tilde{a}_{n+1}}{n}\right|^{2}<\infty
\end{equation*}
and
\begin{equation*}\label{f3}
\left\|\sum\limits_{n\in B^{+}}\tilde{b}_{n}\mathcal{O}(\frac{1}{n})
\tilde{v}(x,\lambda_{n})\right\|_{L^{2}(0,1)}\leq C' \sum_{n\in
B^{+}}|\frac{\tilde{b}_n}{n}|^{2}< \infty.\end{equation*} In a
similar way as in \eqref{bg1} and \eqref{bg2}, by \eqref{bg} and
\eqref{i2}, yields
\begin{align*}
&\left\|\sum\limits_{n\in A\setminus\Lambda^{*}}\tilde{a}_{n}
\xi_{n}\sqrt{\alpha_{n}}
\tilde{v}_{\lambda}(x,\alpha_{n})(1,i\sqrt{\mu_{n-1}})\right\|_{H^{+}}
\leq C_{1}\sum\limits_{n\in A\setminus
\Lambda^{*}}\big|\tilde{a}_{n} \xi_{n}\big|^{2} <\infty,
\end{align*}
\begin{align*}
&\left\|\sum\limits_{n\in A}\tilde{a}_{n+1}
\mathcal{O}(\frac{1}{n})\sqrt{\alpha'_{n}}
\tilde{v}_{\lambda}(x,\alpha'_{n})(1,i\sqrt{\mu_{n-1}})\right\|_{H^{+}}
\leq C_{2}\sum\limits_{n\in A
}\big|\frac{\tilde{a}_{n+1}}{n}\big|^{2} <\infty
\end{align*}
and
\begin{align*}
&\left\|\sum\limits_{n\in B^{+}}\tilde{b}_{n}
\mathcal{O}(\frac{1}{n})\sqrt{\beta_{n}}
\tilde{v}_{\lambda}(x,\beta_{n})(1,i\sqrt{\mu_{n-1}})\right\|_{H^{+}}
\leq C_{3} \sum_{n\in B^{+}}|\frac{\tilde{b}_{n}}{n}|^{2}< \infty.
\end{align*} As $U^{*}_{\big|(0,1)}\in \mathcal{V}_{2}\times
L^{2}(0,1)$, then from the above and \eqref{f1}, we have
\begin{align}
&\Big\|\sum\limits_{n\in \Lambda^{*}}\tilde{a}_{n}
\overline{\tilde{v}}(x,\lambda_{n})+\sum\limits_{n\in
\mathbb{Z}^{*}\setminus ( \Lambda^{*}\cup
B)}\tilde{a}_{n}\overline{\tilde{v}}(x,\mu_{n-1}) +\sum\limits_{n\in
B^{+}}\tilde{b}_{n}\overline{\tilde{v}}(x,\mu_{n-1})\Big\|_{H^{+}}<
\infty.\label{hk}
\end{align}
Recall that $\lambda_{n}=\mu_{n}$ for $n\in \Lambda^{*}$ and
$\{\overline{\tilde{v}}(x,\mu_{n})~\hbox{such
that}~\tilde{v}(0,\mu_{n})=0\}_{n\in \mathbb{Z}^{*}}$ is an
orthogonal subset in $\mathcal{V}_{2}\times L^{2}(0,1)$ with
$\|\overline{\tilde{v}}(x,\mu_{n})\|_{H^{+}}=\mathcal{O}(1)$. Under
the above assumption, since $\mu_{n-1}$ is an eigenvalue of Problem
\eqref{g2}, then in view of \eqref{hk}, we obtain
\begin{equation*}
\sum\limits_{n\in \Lambda^{*}}|\tilde{a}_{n}|^{2}
<\infty,~\sum\limits_{ n\in A\setminus
\Lambda^{*}}|\tilde{a}_{n}|^{2}<\infty,~\sum\limits_{ n\in
A}|\tilde{a}_{n+1}|^{2}<\infty~\hbox{and}~\sum\limits_{n\in
B^{+}}|\tilde{b}_{n}|^{2} <\infty.
\end{equation*}
This implies that $(c_{n})_{n\in A}\in \ell^{2}$ and
$(\tilde{b}_{n})_{n\in B^{+}}\in \ell^{2}$. The proof is complete
now.
\end{proof}
%\begin{rem}
%The asymptotic estimates \eqref{i2} of the Fourier coefficients
%$a_{n}$ are preliminary. More precise asymptotic behaviors of
%$a_{n}$ for $n\in \mathbb{Z}^{*}$ are given at the end of the proof
%of Theorem \ref{xx}.
%\end{rem}
\section{Proof Of Theorem \ref{hj}}
In this section we prove the main controllability result for System
\eqref{a}-\eqref{b1}. Notice that, in view of the fact that
\eqref{a}-\eqref{b1} is linear and reversible in time, this system
is exactly controllable if and only if the system is null
controllable. Using HUM \cite{J.L2}, the controllability of System
\eqref{a}-\eqref{b1} can be reduced to the obtention of an
observability
inequality \eqref{obs} of the uncontrolled System \eqref{a}-\eqref{22}. %%Then the
%%following holds:
%\begin{lemm}\label{lemm1}
%If $T>2 \gamma$, there exist two constants $C_{1}, C_{2}>0$ such
%that
%\begin{equation}\label{ob2}
%C_{1}\|U^{0}\|^{2}_{\mathcal{Y}}\leq \int_{0}^{1}|v_{x}(1,t)|^{2}dt
%\end{equation}
%and
%\begin{equation}\label{ob1}
%\int_{0}^{1}|v_{x}(1,t)|^{2}dt \leq C_{2}
%\|U^{0}\|^{2}_{\mathcal{Y}},
%\end{equation}
%for all $U^{0}=((u^{0},v^{0},z^{0}),(u^{1},v^{1},z^{1}))\in
%\mathcal{Y}$ with $\|.\|_{\mathcal{Y}}$ is defined in \eqref{y}.
%\end{lemm}
As a direct application of Lemma \ref{cor2}, for $T> 2\gamma$ there
exists a positive constant $C$ such that
\begin{align*}
&\int_{0}^{1}| v_{x}(1,t)|^{2}dt=\int_{0}^{1}\sum_{n\in
\mathbb{Z}^{*}}|a_{n}e^{i\sqrt{\lambda_{n}}t}\phi_{n}'(1)|^{2}dt\\
&\geq C \sum_{n \in A} (\sqrt{\la_{n+1}}-\sqrt{\la_{n}})
\left(|a_{n}\phi_{n}'(1)|^{2}+|a_{n+1}\phi_{n+1}'(1)|^{2}\right)\\
&+C\sum_{n \in
A}|a_{n+1}\phi_{n+1}'(1)+a_{n}\phi_{n}'(1)|^{2}+\sum_{n \in
B}|b_{n}\phi_{n}'(1)|^{2}.
\end{align*}
By setting $\tilde{a}_{n}=a_{n}\phi_{n}'(1)$, we obtain
\begin{align*}
&\int_{0}^{1}|v_{x}(1,t)|^{2}dt\geq C\left(\sum_{n \in A}
\delta_{n}^{2}\left(|\tilde{a}_{n}|^{2}+|\tilde{a}_{n+1}|^{2}\right)
+|\tilde{a}_{n+1}+\tilde{a}_{n}|^{2}\right)\\
&+C\sum_{n \in B}|\tilde{b}_{n}|^{2}=C\|U^{0}\|^{2}_{\mathcal{Y}},
\end{align*}
with $U^{0}=((u^{0},v^{0},z^{0}),(u^{1},v^{1},z^{1}))$ and the
asymmetric space $\mathcal{Y}$ is defined by \eqref{y} and
characterized by Theorem \ref{xx}. The second inequality of
\eqref{obs} can be proved in a similar way. This ends up the proof
of Theorem \ref{hj}.
%\begin{equation*}
%\int_{0}^{1}|v_{x}(1,t)|^{2}dt\leq C'\|U^{0}\|^{2}_{\mathcal{Y}}.
%\end{equation*}
%\begin{align*}
%\int_{0}^{1}|v_{x}(1,t)|^{2}dt %\geq C''\left( \sum_{n \in \Lambda}
%%\left[(|a_{n}|^{2}+|a_{n+1}\delta_{n}|^{2})\delta_{n}^{2}+|\delta_{n}a_{n+1}+a_{n}|^{2}\right]\right)\\
%%&+C''\left( \sum_{n \in A\setminus \Lambda}
%%\left[(|a_{n}|^{2}+|a_{n+1}|^{2})\delta_{n}^{4}+\delta_{n}^{2}|a_{n+1}+a_{n}|^{2}\right]
%%+\sum_{n \in B^{+}}|b_{n}|^{2}+\sum_{n
%%\in B^{-}}\frac{|b_{n}|^{2}}{n^{2}}\right)\\\\
%\geq C'' \|U^{0}\|^{2}_{\mathcal{Y}}.
%\end{align*}
%A similar computation shows that \eqref{ob1} holds too.
\begin{rem}\label{rem} Notice that if the coefficients
$\rho_{i},~\sigma_{i},~q_{i}$ $i=1,2$, are symmetric functions, then
the eigenvalues $\mu_{n}^{-}$ and $\mu_{n}^{+}$ of Problems
\eqref{g1} and \eqref{g2} coincide. Moreover, we have
$\lambda_{\varphi(n)}=\lambda_{2n}=\mu_{n}^{+}$, $A=\{2n,~n\in
\mathbb{Z}^{*}\}$ and $B=\emptyset$. Hence, we can use Ullrich's
theorem \cite{D.U}, as in \cite{CAS, E.Z}, to prove the exact
controllability.
\end{rem}
\begin{rem} In the case of Neumann boundary
control, (i.e., $v_{x}(1,t)=p(t)$ and $u(-1,t)=0$), we show that the
associated spectral gap $(\sqrt{\lambda_{n+1}}-\sqrt{\lambda_{n}})$
is uniformly positive if and only if $\frac{\gamma_{1}}{\gamma}$ is
a rational number with an irreductible fraction such that
$\frac{\gamma_{1}}{\gamma}=\frac{2p-1}{q}$, where
$\gamma_{1}=\int_{-1}^{0}(\frac{\rho(x)}{\sigma(x)})^{\frac{1}{2}}dx$,
$\gamma=\int_{-1}^{0}(\frac{\rho_1(x)}{\sigma_1(x)})^{\frac{1}{2}}dx+
\int_{0}^{1}(\frac{\rho_2(x)}{\sigma_2(x)})^{\frac{1}{2}}dx$ and
$p,q\in \mathbb{N}^{*}$. Under this condition, we prove that this
system is exactly controllable in an asymmetric space for a control
time $T> 2\gamma$. This will be the object of a subsequent paper.
\end{rem}
\begin{rem} Note that the method used in this paper
 could be extended to prove the exact boundary controllability of hinged two Euler-Bernoulli
beams connected by a point mass with variable coefficients. The case
of constant coefficients was studied in \cite{Beam2}. In that paper,
it was shown that the associated spectral gap is uniformly positive,
so that by use of Haraux's theorem, it was proved the exact
controllability. It could be conjectured that the results obtained
in \cite{Beam2} remain true in the case of variable coefficients:
all the eigenvalues are simple, the associated spectral gap is
uniformly positive, the first derivative of all the eigenfunctions
do not vanish at the extreme $x=1$.
\end{rem}

%\begin{align}
%\mathcal{Y}&=\big\{U=\sum_{n\in \mathbb{N}^{*}} a_{n}
%\overline{\phi}_{n} \in X_{-\frac{1}{2}} \times
%X_{-1}:~\|U\|^{2}_{\mathcal{Y}}=\sum_{n\in
%\Lambda}\delta_{n}^{2}\left(|a_{n}|^{2}+\frac{|a_{n+1}|^{2}}{n^{2}}\right)+|a_{n}+\frac{a_{n+1}}{n}|^{2}\nonumber\\
%&+\sum_{n\in A \setminus
%\Lambda}\delta_{n}^{2}\left(\frac{|a_{n}|^{2}}{n^{2}}+\frac{|a_{n+1}|^{2}}{n^{2}}\right)+\frac{1}{n^{2}}|a_{n}+a_{n+1}|^{2}
%+\sum_{n\in B^{+}}|b_{n}|^{2}+\sum_{n \in
%B^{-}}|\frac{b_{n}}{n}|^{2}< \infty \big\},\label{y}
%\end{align}

\end{document}